\documentclass[reqno,12pt]{amsart}  

 \usepackage{gensymb}
 
 \usepackage{amsmath}
\usepackage{amssymb} 
\usepackage{mathtools}   
\usepackage{graphicx}
\usepackage{graphics}  
\usepackage{color}
\usepackage{verbatim} 
\usepackage[normalem]{ulem} 
\usepackage[mathscr]{eucal}
\usepackage{upgreek}
\usepackage{enumerate} 
 \usepackage[english]{babel} 
 \usepackage{float}
\floatstyle{boxed}
\restylefloat{figure}

\usepackage[titletoc]{appendix}
\numberwithin{equation}{section}

\usepackage[refpage,noprefix]{nomencl}
\usepackage{nomencl} 

\makenomenclature

\newtheorem{theorem}{Theorem}[section] 
\newtheorem{lemma}[theorem]{Lemma} 
\newtheorem{proposition}[theorem] {Proposition} 
 
\newtheorem{remark}[theorem]  {Remark} 
\newtheorem{definition}[theorem] {Definition}

\theoremstyle{definition}

\DeclareMathAlphabet{\mathpzc}{OT1}{pzc}{m}{it}

\newcommand{\bP} {\boldsymbol{P}}

\newcommand{\bO}{\boldsymbol{\Omega}}

\newcommand{\beeta}{\boldsymbol{\eta}}

\newcommand{\bo}{\boldsymbol{\omega}}
\newcommand{\btau}{\boldsymbol{\tau}}

\newcommand{\bzeta}{\boldsymbol{\zeta}}

\newcommand{\bDel}{\boldsymbol{\Del}}

\renewcommand{\L} {\Lambda} 
\renewcommand{\O} {\Omega}

\def\d{\delta} 
 
\newcommand{\eps}{\varepsilon}

\def\l{\lambda} 
 
\def\om{\omega} 
 
\newfam\Bbbfam 
\font\tenBbb=msbm10 
\font\sevenBbb=msbm7 
\font\fiveBbb=msbm5 
\textfont\Bbbfam=\tenBbb 
\scriptfont\Bbbfam=\sevenBbb 
\scriptscriptfont\Bbbfam=\fiveBbb

\newcommand{\R}     {\mathbb{R}} 
\newcommand{\Z}     {\mathbb{Z}} 
\newcommand{\N}     {\mathbb{N}} 
\renewcommand{\P}   {\mathbb{P}} 
 
\newcommand{\E}     {\mathbb{E}}

\def\1{{\mathchoice {1\mskip-4mu\mathrm l}      
{1\mskip-4mu\mathrm l} 
{1\mskip-4.5mu\mathrm l} {1\mskip-5mu\mathrm l}}} 
\newcommand{\ssup}[1] {{\scriptscriptstyle{({#1}})}} 
\def\comment#1{} 
\newtheoremstyle{thm}{2ex}{2ex}{\itshape\rmfamily}{} 
{\bfseries\rmfamily}{}{1.7ex}{} 
 
\newtheoremstyle{rem}{1.3ex}{1.3ex}{\rmfamily}{} 
{\itshape\rmfamily}{}{1.5ex}{} 
 
\newenvironment{proofsect}[1] 
{\vskip0.1cm\noindent{\scshape #1.}\hskip0.5cm}

\newcommand{\Fcal}   {{\mathcal F }}

\newcommand{\Mcal}   {{\mathcal M }}

\newcommand{\Tscr} {\mathscr{T}}

\newcommand{\Lscr} {\mathscr{L}}
\newcommand{\Zscr}{\mathscr{Z}}

 \newcommand{\ex}{{\rm e}} 
 
\renewcommand{\d}{{\rm d}}

\newcommand{\Leb}{{\rm Leb}}

\newcommand{\Del}{{\operatorname{\sf Del}}}

\newcommand{\Exp}{\mathscr{E}\kern-0.2mm{\operatorname{xp}}}
\newcommand{\Log}{\mathscr{L}\kern-0.2mm{\operatorname{og}}}

\newcommand{\heap}[2]{\genfrac{}{}{0pt}{}{#1}{#2}} 
\newcommand{\pr}{{\operatorname {pr}}}
\renewcommand{\emptyset} {\varnothing}

\makeindex

\begin{document}

\title[\hfill Phase transitions\hfill]
{Phase transitions in Delaunay Potts models}


\author{Stefan Adams and  Michael Eyers}
\address{Mathematics Institute, University of Warwick, Coventry CV4 7AL, United Kingdom}
\email{S.Adams@warwick.ac.uk}

\thanks{}
  

 
\keywords{Delaunay tessellation, Gibbs measures, Random cluster measures, percolation, phase transitions, coarse graining, multi-body interaction}  


\begin{abstract}

We establish phase transitions for certain classes of continuum Delaunay multi-type particle systems (continuum Potts models) with infinite range repulsive interaction between particles of different type. In one class of the Delaunay Potts models studied the repulsive interaction is a triangle (multi-body) interaction  whereas  in the second class the interaction is between pairs (edges) of the Delaunay graph. The result for the edge model is an extension of  finite range results in \cite{BBD04} for the Delaunay graph and in \cite{GH96} for continuum Potts models to an infinite range repulsion  decaying with the edge length. This is a proof of an old conjecture of Lebowitz and Lieb. The repulsive triangle  interactions have infinite range as well and depend on the underlying geometry and thus are a first step towards studying phase transitions for geometry-dependent multi-body systems.
Our approach involves a Delaunay random-cluster representation analogous to the Fortuin-Kasteleyn representation of the Potts model. The phase transitions manifest themselves in the percolation of the corresponding random-cluster model. Our proofs rely on recent studies \cite{DDG12} of Gibbs measures for geometry-dependent interactions.
\end{abstract} 
\maketitle


\section{Introduction}
Although the study of phase transitions is one of the main subjects of mathematical statistical mechanics, examples of models exhibiting phase transition are mainly restricted  to lattice systems. In the continuous setting results are much harder to obtain, e.g., the proof of a liquid-vapor phase transition in \cite{LMP99}, or the spontaneous breaking of rotational symmetry in two dimensions for a Delaunay hard-equilaterality like interaction \cite{MR09}. These phase transitions manifest themselves in  breaking of a continuous symmetry. There is another  specific model for which a phase transition is known to occur: the model of Widom and Rowlinson \cite{WR70}. This is a multi-type particle system in $ \R^d, d\ge 2 $, with hard-core exclusion between particles of different type, and no interaction between particles of the same type. The phase transition in this model was stablished by Ruelle \cite{Rue71}. Lebowitz and Lieb \cite{LL72} extended his result by replacing the hard-core exclusion by a soft-core repulsion between unlike particles. Finally, phase transition results for a general class of continuum Potts models in $ \R^d, d\ge 2 $, have been derived in \cite{GH96}. The repulsive interaction between particles of different type in \cite{GH96} is of finite range and a type-independent background potential has been added. The phase transitions for large activities in all these systems reveal themselves in breaking of the symmetry in the  type-distribution. The results in \cite{GH96} are  based an a random-cluster representation analogous to the Fortuin-Kasteleyn representation of  lattice Potts models, see \cite{GHM}. In \cite{BBD04} the soft repulsion in \cite{GH96} between unlike particles has been replaced by another kind of soft repulsion based on the structure of some graph. More precisely, the finite range repulsion between particles of different type acts only on a nearest neighbour subgraph of the Delaunay graph in $ \R^2$. 

In this paper we establish the existence of phase transitions for two different classes of continuum Delaunay Potts models in $ \R^2 $ including the Widom-Rowlinson model. The repulsive interaction between unlike particles in all models is of infinite range, and it depends in one class on the geometry of the Delaunay triangle (minimal angle) and in another class on the Delaunay edge when the interaction decays with the length of the edge. All our models come with a hard-core background interaction which is independent of the type of the particles. We believe that one can extend our results by removing this hard-core constraint, see Remark~\ref{softcore} below. Our results extend \cite{BBD04} and \cite{GH96} in two ways by having an infinite range repulsion and a geometry-dependency of the type interaction. Gibbs models on Delaunay structures have been studied in \cite{BBD99,BBD02,BBD04,Der08,DDG12,DG09,DL11}, and our results rely on the existence of Gibbs measures for geometry-dependent interactions established  in \cite{DDG12}. Our approach is based on a Delaunay random-cluster representation similar to \cite{BBD04} and \cite{GH96}, the difference being that we replace edge percolation by an adaptation of lattice  hyperedge percolation \cite{G94} to our continuum setting. A phase transition for our Delaunay Potts models follows if we can show that the corresponding percolation process contains an infinite cluster. A similar program was carried out by Chayes et al. in \cite{CCK} for the hard-core Widom-Rowlinson model. In that case, the existence of infinite clusters follows from a stochastic comparison with the Poisson Boolean model of continuum percolation, while our framework uses a coarse graining method to derive a stochastic comparison with site percolation on $ \Z^2 $. This idea actually goes  back to \cite{Haggstrom00} in which percolation for the Poisson Voronoi model has been established. We note that our random-cluster representation requires the symmetry of the type  interaction. In the non-symmetric Widom-Rowlinson models, the existence of a phase transition has been established by Bricmont et al. \cite{BKL}, and recently by Suhov et al. \cite{MSS}. 

We conclude with some remarks on the particular features models defined on the Delaunay hypergraph structure show. The most simple Delaunay Potts model would be one with no more interaction than constant Delaunay edge interaction as then the percolation would be independent of the activity parameter. Our triangle model  is a  step towards understanding this system.  We draw attention to some differences between geometric models on the Delaunay hypergraph structure and that of classical models such as the Widom-Rowlinson model and its soft-core variant of Lebowitz and Lieb \cite{LL72}. The first is that edges and triangles in the Delaunay hypergraph are each proportional in number to the number of particles in the configuration. However, in the case of the complete hypergraph the number of edges is proportional to the number of particles squared and the number of triangles is proportional to the number of particles cubed. Secondly, in the complete graph of the classical models, the neighbourhood of a given point depends only on the distance between points and so the number of neighbours increases with the activity parameter $ z $ of the underlying point process. This means that the system will become strongly connected for high values of $z$. This is not the case for the Delaunay hypergraphs which exhibit a self-similar property. Essentially, as the activity parameter $z$ increases, the expected number of neighbours to a given point in the Delaunay hypergraph remains the same, see \cite{M94}. Therefore, in order to keep a strong connectivity in our geometric models on Delaunay hypergraphs, we use a type interaction between particles of a hyperedge with a non-constant mark. Finally, and perhaps most importantly, is the question of additivity. Namely, suppose we have an existing particle configuration $\om$ and we want to add a new particle $x$ to it. In the case of classical many-body interactions, this addition will introduce new interactions that occur between $x$ and the existing configuration $\om$. However, the interactions between particles of $\om $ remain unaffected, and so classical many-body interactions are additive. On the other hand, in the Delaunay framework, the introduction of a new particle to an existing configuration not only creates new edges and triangles, but destroys some too. The Delaunay interactions are therefore not additive, and for this reason, attractive and repulsive interactions are indistinct. In the case of a hard exclusion interaction, we arrive at the possibility that a configuration $ \om $ is excluded, but for some $x$, $ \om\cup x $ is not. This is called the non-hereditary property \cite{DG09}, which seems to rule out using techniques such as stochastic comparisons of point processes \cite{GK97}.

\section{Results}
\subsection{Setup}
We consider configurations of points in $ \R^2 $ with internal degrees of freedom, or marks. Let $ E=\{1,\ldots,q\}, q\in\N, q\ge 2 $, be the finite set of different marks. That is, each marked point is represented by a position $ x\in\R^2 $ and a mark $ \sigma(x)\in E $, and each marked configuration $ \bo $ is a countable subset of $ \R^2\times E $ having a locally finite projection onto $ \R^2 $. We denote by $ \bO $ the set of all marked configurations with locally finite projection onto $ \R^2 $. We will sometimes identify $ \bo $ with a vector $ \bo=(\omega^{\ssup{1}},\ldots,\om^{\ssup{q}}) $ of pairwise disjoint locally finite sets $ \om^{\ssup{1}},\ldots,\om^{\ssup{q}} $ in $ \R^2 $ (we write $ \O $ for the set of all locally finite configurations in $ \R^2 $). Any $ \bo $ is uniquely determined by the pair $ (\om,\sigma)$, where $ \omega=\cup_{i=1}^q\omega^{\ssup{i}} $ is the set of all occupied positions, and where the mark   function $ \sigma\colon \om\to E $ is defined by $ \sigma(x)=i $ if $ x\in\om^{\ssup{i}}, i\in E $. For each measurable set $B$  in $ \R^2\times E $ the counting variable $ N(B)\colon\bo\to\bo(B) $ on $ \bO $ gives the number of marked particles such that  the pair (position, mark) belongs to $B$. We equip the space $ \bO $ with the $\sigma$-algebra $\boldsymbol{\Fcal} $ generated by the counting variables $N(B) $ and the space  $ \O $ of locally finite configurations with the $\sigma$-algebra $ \Fcal $ generated by the counting variables $ N_\Delta=\#\{\om\cap\Delta\} $ for $ \Delta\Subset\R^2 $ where we write $ \Delta\Subset\R^2 $ for any bounded $ \Delta\subset\R^2 $. As usual, we take as reference measure on $ (\bO,\boldsymbol{\Fcal}) $ the marked Poisson point process $ \boldsymbol{\Pi}^z $ with intensity measure $ z\Leb\otimes\mu_{\sf u} $  where $ z>0 $ is an arbitrary activity, $ \Leb $ is the Lebesgue measure in $ \R^2 $, and $ \mu_{\sf u} $ is the uniform probability measure on $ E $.


We let $ \bO_f \subset \bO$ (resp. $ \O_f\subset\O $) denote the set of all finite configurations.  For each $ \L\subset\R^2 $ we write $ \bO_\L=\{\bo\in\bO\colon \bo\subset\L\times E\} $  for the set of configurations in $ \L $, $  \pr_\L\colon\bo\to\bo_\L:=\bo\cap \L\times E $ for the projection from $ \bO $ to $ \bO_\L $ (similarly for unmarked configurations), $\boldsymbol{\Fcal}_\L^\prime=\boldsymbol{\Fcal}|_{\bO_\L} $ for the trace $\sigma$-algebra of $ \boldsymbol{\Fcal} $ on $ \bO_\L $, and  $ \boldsymbol{\Fcal}_\L=\pr_\L^{-1}\boldsymbol{\Fcal}_\L^\prime\subset\boldsymbol{\Fcal} $ for the $\sigma$-algebra of all events that happen in $ \L $ only. The reference measure on $ (\bO_\L,\boldsymbol{\Fcal}_\L^\prime) $ is $ \boldsymbol{\Pi}_\L^z:=\boldsymbol{\Pi}^z\circ\pr_\L^{-1} $. In a similar way we define the corresponding objects for unmarked configurations, $ \Pi^z, \Pi^z_\L, \O_\L, \pr_\L, \Fcal_\L^\prime $, and $ \Fcal_\L $. Finally, let $ \varTheta=(\vartheta_x)_{x\in\R^2} $ be the shift group, where $ \vartheta_x\colon\bO\to\bO $ is the translation of the  spatial component by the vector $ -x\in\R^2 $. Note that by definition, $ N_\Delta(\vartheta_x\bo)=N_{\Delta+x}(\bo) $ for all $ \Delta\subset\R^2 $. \medskip

The interaction between the  points in all the models to be studied  depend on the geometry of their location. We  describe this in terms of hypergraph structures. A hypergraph structure is a measurable set of $ \O_f\times\O $ (resp. $ \bO_f\times\bO $). We outline the definitions for the unmarked configurations first with obvious adaptations to the case of marked point configurations.  The set $ \Del $ of Delaunay hyperedges consist of all pairs $ (\eta,\om)$ with $ \eta\subset\om $ for which there exits an open ball $ B(\eta,\om) $ with $ \partial B(\eta,\om)\cap\om=\eta $ that contains no points of $ \om $. For $ m=1,2,3  $, we write $\Del_m=\{(\eta,\om)\in\Del\colon \#\eta=m\} $ for the set of Delaunay simplices with $ m $ vertices. Given a configuration $\om $  the set of all Delaunay hyperedges $ \eta\subset\om $ with $ \#\eta=m $ is denoted by $\Del_m(\om) $. It is possible that $ \eta\in\Del(\om) $ consists of four or more points on a sphere with no points inside. In fact, for this not to happen, we must consider configurations in general   position as in \cite{M94}. More precisely, this means that no four points lie on the boundary of a circle and every half-plane contains at least one point. Fortunately, this occurs with probability one for our Poisson reference measures, and in fact, for any stationary point process. Note that the open ball $ B(\eta,\om) $ is only uniquely determined when $ \#\eta=3 $ and $ \eta $ is affinely independent. Henceforth, for each configuration $ \om $ we have an associated Delaunay triangulation 
\begin{equation}\label{triangulation}
\{\tau\subset \om\colon \#\tau=3,B(\tau,\om)\cap\om=\emptyset\}
\end{equation}
of the plane, where $ B(\tau,\om) $ is the unique open ball with $ \tau\subset\partial B(\tau,\om) $. The set in \eqref{triangulation} is uniquely determined and defines a triangulation of the convex hull of $ \om $ whenever $ \om $ is in general position (\cite{M94}). In a similar way one can define the marked Delaunay hyperedges, i.e., $ \boldsymbol{\Del} $ and $ \boldsymbol{\Del}_m(\bo) $ as measurable sets in $ \bO_f\times\bO $ where the Delaunay property refers to the spatial component only.

\medskip

Given a configuration $ \om \in\Omega $ (or $ \bo$) we write $ \O_{\L,\om}=\{\zeta\in\O\colon \zeta\setminus\L=\om\} $ (resp. $ \bO_{\L,\bo} $) for the set of configurations which equal $ \om $ off $ \L $.  For any triangle $ \tau\in\Del_3 $ we denote the minimal angle by $ \alpha(\tau)\in(0,\pi/3] $ and for any edge $ \eta\in\Del_2 $ we denote its  length by $ \ell(\eta) $. We write $ \L\Subset\R^2 $ in the following when $ \L $ is non-empty and finite, i.e., $ 0<|\L|<\infty $.
The interaction is given by the following Hamiltonian in $ \L $ with boundary condition $ \bo\in\bO $, $m=2,3$, 
\begin{equation}\label{Hamiltonian}
H_{\L,\bo}(\bzeta):=\sum_{\eta\in\Del_{2,\L}(\zeta)}\psi(x-y)+\sum_{\heap{\btau\in\bDel_m(\bzeta),}{\tau\in\Del_{m,\L}(\zeta)}}\phi(\alpha(\tau))(1-\delta_\sigma(\btau)),\quad\bzeta\in\O_{\L,\bo},
\end{equation}
where $ \Del_{2,\L}(\zeta):=\{\eta\in\Del_2(\zeta)\colon \,\exists\,  \tau\in\Del_3(\zeta), \eta\subset\tau,\partial B(\tau,\zeta)\cap\L\not=\emptyset\} $ and $ \Del_{3,\L}(\zeta):=\{\tau\in \Del_3(\zeta)\colon \partial B(\tau,\zeta)\cap \L\not=\emptyset\} $.
Here $ \psi\colon\R^2 \to\R\cup\{\infty\} $ is an even measurable function and  $ \phi $ is a measurable function of the minimal angle $ \alpha(\tau) $ when $ m=3 $ or a measurable function  of the length $ \ell(\eta) $ of an edge when $m=2 $ (see below) and 
$$ 
\delta_\sigma(\btau)=\begin{cases} 1 &, \mbox{ if } \sigma_{\btau}(x)=\sigma_{\btau}(y) \mbox{ for any pair } \{x,y\}\subset\tau,\\
0&, \mbox{ otherwise}.
\end{cases}
$$ 
The second term on the right hand side of \eqref{Hamiltonian} is the type interaction and thus the most important contribution. It describes a repulsion between particles having a different mark. The first term corresponds to a 'colour-blind' background  interaction. Specifically, we assume the following for our models.

\smallskip

\noindent\textbf{Assumptions:}

 \begin{enumerate}
 \item[{(A1)}] 
 \begin{itemize} 
 \item[(i)] \textbf{Triangle-model I}: $m=3$ and $ \phi\colon(0,\pi/3]\to\R_+ $ and there are $ \alpha_0\in(0,\pi/3) $ and $ \beta> 0 $ such that
 $$
 \phi(\theta)=\begin{cases}  0&, \mbox{ if } \theta< \alpha_0,\\ \beta&, \mbox{ if } \theta\ge \alpha_0. \end{cases}
 $$ 
 
 \item[(ii)] \textbf{Triangle-model II:} $m=3 $ and $ \phi\colon(0,\pi/3]\to\R_+ $ with 
 $$ 
 \phi(\theta)=\log(1+\beta\theta^3),\qquad \beta \ge 0.
 $$
 
 \item[(iii)] \textbf{Edge-model:} Let $ \delta_0>0 $ be some parameter to be specified in (A2), $m=2$ and $\phi\colon[0,\infty]\to\R_+\cup\{\infty\} $ with 
 $$
 \phi(\ell)=\log(1+\beta(\delta_0/\ell)^3),\qquad \beta\ge 0.
 $$
 
 \end{itemize}
 \item[{(A2)}] $ \psi $ is a hard-core potential, that is, there is a range $ \delta_0>0 $ such that $ \psi(x-y)=0 $ whenever $ |x-y|\ge \delta_0 $ and $ \Psi(x-y)=+\infty $ otherwise.

  \end{enumerate}
 Following \cite{DDG12} we say a  configuration $ \bo\in\bO $   (or $\om\in\O$) is \textbf{admissible} for $ \L\Subset\R^2 $ and activity $ z$ if $ H_{\L,\bo} $ is $\boldsymbol{\Pi}^z$-almost surely well-defined and $ 0<Z_\L(\bo)<\infty $, where the  partition function is defined as
  $$
  Z_\L(\bo)=\int_{\bO_{\L,\bo}}\,{\rm e}^{-H_{\L,\bo}(\bzeta)}\,\Pi^z(\d\zeta^{\ssup{1}})\cdots\Pi^z(\d\zeta^{\ssup{q}}).
  $$
 We denote the set of admissible configurations by $ \bO_\L^* $. The Gibbs distribution for $ \psi, \phi $, and $ z> 0 $ in $\L $ with admissible boundary condition $ \bo $ is defined as
 \begin{equation}\label{Gibbsdist}
 \gamma_{\L,\bo}(A)=\frac{1}{Z_\L(\bo)}\int_{\bO_{\L,\bo}}\,\1_A(\bzeta\cup\bo){\rm e}^{-H_{\L,\bo}(\bzeta)}\,\Pi_\L^z(\d\bzeta),\quad A\in\boldsymbol{\Fcal}.
 \end{equation}
  It is evident from \eqref{Gibbsdist} that, for fixed $ \zeta\in\O_\L $, the conditional distribution of the marks of $ \bzeta= (\zeta^{\ssup{1}},\ldots,\zeta^{\ssup{q}})$  relative to $ \gamma_{\L,\bo} $ is that of a discrete Potts model on $ \zeta $ embedded in the Delaunay triangulation with position-dependent interaction between the marks. This justifies calling our models  \textit{Delaunay Potts models}.
  
  \begin{definition}
  A probability measure $ \mu $ on $ \bO $ is called a Gibbs measure for the Delaunay  Potts model with activity $ z>0 $ and interaction potentials $ \psi$ and $ \phi $  if $ \mu(\bO^*_\L)=1 $ and
 \begin{equation}\label{DLR}
   \E_\mu[f]=\int_{\bO^*_\L}\,\frac{1}{Z_\L(\bo)}\int_{\bO_{\L,\bo}}f(\bzeta\cup\bo){\rm e}^{-H_{\L,\bo}(\bzeta)}\,\boldsymbol{\Pi}^z_\L(\d\bzeta)\mu(\d\bo)
 \end{equation} for every $ \L\Subset\R^2 $ and every measurable function $ f$.
  \end{definition}
  The equations in \eqref{DLR} are the DLR equations (after Dobrushin, Lanford, and Ruelle). They express that the Gibbs distribution in \eqref{Gibbsdist} is a version of the conditional probability $ \mu(A|\boldsymbol{\Fcal}_{\L^{\rm c}})(\bo) $. The measurability of all objects is established in \cite{E14,DDG12}.

\subsection{Results and remarks}   
 We outline the results for our triangle and edge  models all of which are generic examples of a more general class of coarse grain ready potentials, see Section~\ref{tileperc} for further details.

  \begin{proposition}[\textbf{Existence of Gibbs measures}]\label{Propexist}
  
  There exist at least one Gibbs measure for the following Delaunay Potts  models with hard-core background potential of range $ \delta_0>0 $.  
 \begin{enumerate} 
 
 \item[(a)] 
 \begin{itemize}
 
 \item[(i)] \textbf{Triangle model I:} For any $ z> 0 $ and for any $ \alpha_0\in(0,\pi/3) $ and $ \beta\ge 0 $. 
 \item[(ii)] \textbf{Triangle model II:} For any $ z> 0 $  and $ \beta\ge 0 $.
 \end{itemize}
 
 \medskip

 \item[(b)] \textbf{Edge model:} For any $ z> 0 $ and $ \beta\ge 0 $.
 
  \end{enumerate}
 \end{proposition}

 \begin{remark}[\textbf{Gibbs measures}]
 We obtain the existence of Gibbs measures for large enough activity only. This is mainly due to the techniques used in \cite{DDG12} for the existence of Gibbs measures (\cite{G88}). The restriction originates from the need to obtain an upper bound for the Hamiltonian for some chosen configurations, namely, the so-called pseudo-periodic configurations (see Appendix~\ref{pseudoperiodic} or \cite{DDG12}).  Large activities then ensure that these pseudo-periodic configurations have sufficient mass under the reference process.  Existence of Gibbs measures for related Delaunay models have been obtained in \cite{BBD99, Der08,DG09}. Note that for $q=1$ our models have no marks and Gibbs measures do exist as well (\cite{DDG12}).
 \end{remark}
 
 Gibbs measures for the Delaunay Potts model do exist for large enough activity $z>0 $. A phase transition is said to occur if there exists more than one Gibbs measure for the Delaunay Potts model. The following two theorems show that this happens when the activity $ z $ and the parameter $ \beta $ are sufficiently large. Note that $ \beta $ is a parameter for the type interaction and not the usual inverse temperature. The first theorem extends the results in \cite{BBD04} for finite range type interaction  to the case of infinite range type interaction.
  
  \begin{theorem}[\textbf{Phase transition - Edge model}]\label{THM-mainedge}
  For all $ \delta_0>0 $ there is  $ \beta_0=\beta_0(\delta_0) $ such that for all $ \beta>\beta_0 $ there is  $ z_0=z_0(\beta,\delta_0) $ such that for  all $ z>z_0 $, there exit at least $q$ different  Gibbs measures for the Delaunay edge model. 
  \end{theorem}
  
  \begin{remark}[\textbf{Range of type-interaction}]
 The edge model is an infinite range   type-interaction while keeping a hard-core background potential as in \cite{BBD04}, where it is used on one hand to obtain a bound on the number of connected components in the hypergraph $ \Del_2 $ once a point is inserted into a configuration, and on the other hand, to bound the number of points in a coarse graining scheme. We have replaced the bound on the number of connected components by a bound of the  expectation under the edge drawing measure using the hard-core background potential and an sufficient decay of the type interaction. We believe that one can relax the hard-core assumption, see Remark~\ref{softcore}. \end{remark}

 The following theorem shows that phase transitions occur in our Delaunay triangle models. As outlined in \cite[Remark~3.11]{DG09} a triangle model similar to the edge model in  \cite{BBD04} shows phase transition as well. The triangle interaction in \cite{DG09}  is of finite range as well. Our triangle models have an infinite range type interaction. We obtain phase transitions for sufficiently large activity and parameter $ \beta $ for two models with infinite range  type interactions. The idea for the triangle model I is that tiles (triangles) whose minimal angles exceed a given threshold $ \alpha_0 $ gain in probability when their vertices are of the same type. The energy of a triangle with minimal angle is given by the constant value of the type potential $ 
\phi $. In our second triangle model II,  however, this incentive is not constant as it grows with the minimal angle and it decays with shrinking minimal angle.

  \begin{theorem}[\textbf{Phase transition - Triangle model}]\label{THM-main}
  The following holds for the triangle models.
  \medskip

  \begin{itemize}

  \item[(a)] \textbf{Triangle Model I:} For all $ \delta_0>0 $ and $ \alpha_0 $ sufficiently small there are $ z_0=z_0(\delta,\alpha_0) $ and $ \beta_0=\beta_0(\delta_0,\alpha_0) >0 $ such that for all $ z>z_0 $ and $ \beta>\beta_0 $ there exists at least $q$ different Gibbs measures for the Delaunay triangle model I. 
  
  \item[(b)] \textbf{Triangle Model II:} For all $ \delta_0>0 $ there is $ \beta_0=\beta_0(\delta_0) $ such that for all $ \beta>\beta_0 $ there is  $ z_0=z_0(\beta,\delta_0) $ such that for  all $ z>z_0 $, there exit at least $q$ different  Gibbs measures for the Delaunay triangle model II. 
   \end{itemize}
  \end{theorem}
  
For the triangle model I note that our proof techniques require the minimal angle threshold $ \alpha_0 $ to be sufficiently small, see Lemma~\ref{LemTMI}. In contrast to this, the triangle model II allows the smallest angles to take any value as smaller values experience vanishing preference for type alignment in a triangle.

\begin{remark}[\textbf{Hardcore versus softcore background potential}]\label{softcore}
In all our models we add a hard-core background potential between pairs of particles. This is for technical reasons only, and we believe that it is possible to replace it  by some  triangle potential excluding minimal angles below a given threshold (see \cite{DDG12}), or by the hard-equilaterality model (\cite{DDG12},\cite{DL11}). These triangle background potentials are of infinite range in terms of Delaunay edges (pairs of particles).  They allow to bound easily the number of connected components in the Delaunay random-cluster model once a new point is inserted. On the other hand inserting a point in an empty cell now depends on all the points in neighbouring cells. An adaptation of our  coarse graining method in Section~\ref{tileperc}  and restricting  to certain admissible boundary conditions for the cells allows to estimate the probability of injecting a  single particle.  Another possibility is to replace the hard-core background pair potential by a superstable pair potential (soft-core background) as in \cite{R70},\cite{GH96}, and \cite{BBD02}. This allows to control the number of particles in cells and  to estimate the probability of insertion of single points, see e.g. \cite{GH96}. However, for the final percolation argument it seems an additional  multi scale approach is needed, which is based in our coarse graining technique in Section~\ref{tileperc}. Suppose that in each of the sub cells (paralleloptopes) there is a region with a finite number of points. Then, due to the chosen  partition of $ \R^2$, there will be Delaunay triangles connecting these regions in the different cells. These triangles have then a minimal angle and will thus be open with high probability. Thus all the regions in the different cells are connected with open triangles leaving the questions whether within these regions all triangle are open as well. We think that this can be achieved by some kind a renormalisation approach with a combined bond-site percolation  as in \cite{GH96}. We hope to approach this idea using  Delaunay potentials studied in \cite{DDG12} replacing the standard pair potential in the future (\cite{A15}). 
\end{remark}

\begin{remark}[\textbf{Free energy}]
One may wonder if the phase transitions manifest itself thermodynamically by a non-differentiability ("discontinuity") of the free energy (pressure).  Using the techniques from \cite{Geo94} and  \cite{DG09}, it should be possible to obtain a variational representation of the free energy, see also \cite{ACK11} for free energy representations for marked configurations. Then a discontinuity of the free energy can be established using our results above. For continuum Potts models this has been established in \cite[Remark~4.3]{GH96}. For a class of bounded triangle potential \cite{DG09} shows that the Gibbs measures are minimisers of the free energy. The free energy is the level-3 large deviation rate functional for empirical processes, see \cite{Geo94,ACK11}.
\end{remark}

 \begin{remark}[\textbf{Uniqueness of Gibbs measures}]
 To establish uniqueness of the Gibbs measure in our Delaunay Potts models one can use the Delaunay random-cluster measure $ C_{\L_n,\omega} $, to be defined in \eqref{DRCmeasure} below. In \cite[Theorem~6.10]{GHM} uniqueness is established once the probability of an open connection of the origin to infinity is vanishing for the limiting lattice version of the random-cluster measure, that is, for some set $ \Delta\Subset\R^2$ containing the origin, 
 $$
 \lim_{n\to\infty} C_{\L_n,\om}(\Delta\longleftrightarrow \L_n^{\rm c})=0,
 $$  for a sequence of boxes $ \L_n\Subset \R^2$ with $ \L_n\uparrow\R^2 $ as $ n\to\infty $. One way to achieve this, is to obtain an stochastic domination of the Delaunay  random-cluster measure by the so-called random Delaunay edge model of hard-core particles.  Using \cite{BBD02} we know that the critical probabilities for both, the site and bond percolation on the Delaunay graph, are bounded from below. Extension to our tile (hyperedge) percolation using \cite{G94} can provide a corresponding lower bound as well. Thus, if our parameter $ \beta $ is chosen sufficiently small, then there is no percolation in our Delaunay random-cluster measures and therefore uniqueness of the Gibbs measure.
 \end{remark}

 It goes without saying that the results addressed here are merely a first step towards a closer study of phase transitions with geometry-dependent interactions. The study for Widom-Rowlinson or Potts models with geometry-dependent interaction is by far not complete, one may wish to extend the single tile (edge or triangle) interaction to mutual adjacent Voronoi cell interaction.  The common feature of all these 'ferromagnetic' systems is that phase transitions are due to breaking the symmetry  of the type  distribution. Breaking of continuous symmetries is a much harder business, see e.g. \cite{MR09}, which shows a breaking of rotational symmetry in two dimensions, and which can be seen as a model of oriented particles in $ \R^2 $ with a Delaunay hard-equilaterality interaction. Models studied in \cite{DDG12} may be natural candidates for the existence of a crystallisation transition. A recent ground state study (zero temperature)  \cite{BPT} for multi-body interactions of Wasserstein types shows optimality of the triangular lattice. It seems to be  promising to analyse this model for non-zero temperature in terms of crystallisation transition. We will investigate this further.
 
 Another direction for the class of  'ferromagnetic' models with geometry-dependent interaction is to analyse closer the percolation phenomena, and in particular to study the conformal invariance of crossing probabilities as done for the Voronoi percolation (\cite{BS98,BR06}, and \cite{T14}). 
 
Last but not least it is interesting to study the higher dimensional cases in $ \R^d, d\ge 3 $, as well (see \cite{GH96}). We believe that our results still hold for these cases but there are some technical issues related to the general quadratic position such that we defer that analysis for future study.

 The rest of the paper is organised a follows. In Section~\ref{DelCluster} we define the Delaunay random-cluster measure for tiles (edges/triangles), and in Section~\ref{tileperc} we prove tile percolation for our different  Delaunay random-cluster measures for sufficiently large activities and sufficiently large $ \beta $. In Section~\ref{proofs} we gives details of our proofs.

\section{Delaunay Random Cluster measure}\label{DelCluster}
 For $ \L\Subset\R^2 $  and parameters $ z,\psi,\phi $, and $ \beta $ we define a joint distribution of the Delaunay Potts model and a tile (edge/triangle) process which we call Delaunay random-cluster model. Here tile is used as a general name for any hyperedge of a Delaunay graph,  and  it will be clear from the context if we refer to Delaunay edges or triangles. We follow ideas in \cite{GH96} and will replace edge percolation by tile (i.e., triangle)  percolation.  The basic idea is to introduce random tiles  between points in the plane.   Let 
 $$
 T_{\R^2}=\{\tau=\{x,y,z\}\subset\R^2\colon x\not=y, x\not=z, z\not= y\}
 $$  be the set of all possible tiles of triples of points in $ \R^2 $, likewise, let $T_\L $ be  the set of all tiles  in $ \L $ and $ T_\zeta $ for the  set of tiles in $ \zeta\in\Omega_{\L,\om} $. We identify $ \om $ with $\om^{\ssup{1}} $ and $ \bo=(\om^{\ssup{1}},\emptyset,\ldots,\emptyset) $.  This allows only monochromatic boundary conditions whereas the general version involves the so-called Edwards-Sokal coupling (see \cite{GHM} for lattice Potts models). We restrict ourself to the former case for ease of notation. We write $$ \Tscr =\{T\subset T_{\R^2}\colon T\mbox{ locally finite}\}$$ for the set of all locally finite tile configurations.

 The joint distribution is built from the following three components.

 \noindent The \textit{point distribution}  is  given by the Gibbs distribution $ P_{\L,\om}^{zq} $ in $ \L $ with admissible boundary condition $ \om\in\O_{\L}^*$, interaction $ \psi $, and activity $ zq$, i.e., 
 \begin{equation}\label{blind}
 P_{\L,\omega}^{zq}(\d \zeta)=\frac{1}{Z_\L(\omega)}\exp\Big(-\sum_{\eta\in\Del_2(\zeta)}\psi(\ell(\eta))\Big)\Pi^{zq}_{\L}(\d\zeta),
 \end{equation} with $\zeta\in\Omega_{\L,\om} $ and where $ \ell(\eta)=|x-y| $ is the length of the Delaunay edge $ \eta=\{x,y\} $. Note that for $q=1 $ this measure coincides with \eqref{Gibbsdist}. 
 
 \medskip
 
 \noindent The \textit{type picking mechanism} for a given configuration $ \zeta\in\Omega_{\L,\om} $ is the distribution $ \l_{\zeta,\L} $ of the mark vector $\sigma\in E^\zeta $.  Here $ (\sigma(x))_{x\in\zeta} $ are independent and uniformly distributed random variables on $E$ with $ \sigma(x)=1 $ for all $ x\in\zeta_{\L^{\rm c}}=\om $.  The latter condition ensures that all points outside of $ \L$ carry the given fixed mark.
 \medskip
 
 \noindent The \textit{tile drawing mechanism}. Given a point configuration $ \zeta\in\Omega_{\L,\om} $   we let $ \mu_{\zeta,\L} $ be the distribution of the random tile configuration $ \{\tau\in T_\zeta\colon\upsilon(\tau)=1\} \in\Tscr $ with the hyperedge configuration $ \upsilon\in \{0,1\}^{T_\zeta} $ having
 probability   
 $$
\prod_{\tau\in T_\zeta}p(\tau)^{\upsilon(\tau)}(1-p(\tau))^{1-\upsilon(\tau)}
$$ with
\begin{equation}\label{tiledrawing}
p(\tau):=\P(\upsilon(\tau)=1)=\begin{cases} (1-\ex^{-\phi(\tau)})\1_{\Del_3(\zeta)}(\tau) & \mbox{ if } \tau\in T_{\R^2}\setminus T_{\L^{\rm c}},\\
\1_{\Del_3(\zeta)}(\tau) & \mbox{ if } \tau\in T_{\L^{\rm c}}.
\end{cases} 
\end{equation}
 
The measure $ \mu_{\zeta,\L} $ is a point process on $T_{\R^2} $. Note that $ \zeta\rightarrow\l_{\zeta,\L} $ and $ \zeta\rightarrow\mu_{\zeta,\L} $ are probability kernels (see \cite{E14}). Let the measure
 $$
 P^{zq}_{\L,\bo}(\d\bzeta,\d T)=\frac{1}{Z_\L(\om)}P_{\L,\om}^{zq}(\d\zeta)\l_{\zeta,\L}(\d\bzeta)\mu_{\zeta,\L}(\d T)
 $$ be supported on the set of all $ (\bzeta,T) $ with $ \bzeta\in\bO_{\L,\bo} $ and $ T\subset T_\zeta $. We shall condition on the event that the marks of the points are constant on each connected component in the hypergraph $ (\zeta,T\cap T_\zeta) $. Two distinct vertices $x $ and $y$ in the hypergraph are adjacent to one another if there exists $ \tau\in T_\zeta $ such that $ \{x,y\}\subset\tau $. A path in the hypergraph $ (\zeta,T\cap T_\zeta) $ is an alternating sequence $ v_1,t_1,v_2,t_2,\ldots $ of distinct vertices $ v_i $ and hyperedges (tiles) $t_j $ such that $ \{v_i,v_{i+1}\}\subset t_i $ for all $ i\ge 1 $.  We write
 $$
 A=\{(\bzeta,T)\in\bO\times\Tscr\colon\sum_{\tau\in T}(1-\delta_{\sigma}(\btau))=0\}
 $$ for the set of marked  point configurations such that all vertices of the tiles carry the same mark. The set $ A $ is measurable which one can see from writing the condition in the following way
 $$
 \sum_{\tau=\{x,y,z\}\in T}\;\sum_{i=1}^q\Big(\1_{\zeta^{\ssup{i}}}(x)(1-\1_{\zeta^{\ssup{i}}}(z))+\1_{\zeta^{\ssup{i}}}(y)(1-\1_{\zeta^{\ssup{i}}}(z))+ \1_{\zeta^{\ssup{i}}}(x)(1-\1_{\zeta^{\ssup{i}}}(y))\Big)=0
 $$ and using the fact that $ (\bzeta,x)\mapsto \1_{\zeta^{\ssup{i}}}(x) , i=1,\ldots, q $, are measurable (see \cite[Chapter 2]{GH96}). Furthermore,  $ P^{zq}_{\L,\bo}(A)>0 $, which follows easily observing
 $ P^{zq}_{\L,\bo}(A)\ge P^{zq}_{\L,\om}(\{\widetilde{\om}\})=\ex^{-zq|\L|}/Z_\L(\om) $ where $ \widetilde{\om} $ is the configuration  which  equals $ \om $ outside of $ \L $ and which is empty inside $ \L $.  Henceforth, the random-cluster measure
$$
 \boldsymbol{P}=P^{zq}_{\L,\bo}(\cdot\vert A)
$$
 is well-defined. As in \cite{GH96} we obtain the following two measures from the  random-cluster measure $\bP$, namely if we disregard the tiles we obtain the Delaunay Gibbs distribution $ \gamma_{\L,\bo} $ in \eqref{Gibbsdist} (see \cite{E14}). For the second measure consider the  mapping $ \mbox{sp}\colon(\bzeta,T)\to(\zeta,T) $ from $ \bO\times\Tscr $ onto $ \O\times\Tscr $ where $ \bzeta\mapsto \zeta=\cup_{i=1}^q\zeta^{\ssup{i}} $. For each $ (\zeta,T) $ with $ T\subset T_\zeta $ we let $ K(\zeta,T)$ denote the number of connected components in the hypergraph $ (\zeta,T) $. The Delaunay random-cluster distribution on $ \O\times\Tscr $ is defined by
\begin{equation}\label{DRCmeasure}
  C_{\L,\om}(\d\zeta,\d T)=\frac{1}{\Zscr_\L(\om)}q^{K(\zeta,T)}P^z_{\L,\om}(\d\zeta)\mu_{\zeta,\L}(\d T),
 \end{equation}
 where $ P^z_{\L,\om} $ is \eqref{blind} with activity $ z$ replacing  $zq$ and where
 $$
 \Zscr_\L(\om)=\int_{\O_{\L,\om}}\,\int_{\Tscr}\, q^{K(\zeta,T)}P_{\L,\om}^z(\d\zeta)\mu_{\zeta,\L}(\d T)
 $$ is the normalisation. It is straightforward to  show that $ \bP\circ\mbox{sp}^{-1}=C_{\L,\om} $ (see \cite{E14}). To relate the influence of the boundary condition on the mark of a single point  to the connectivity probabilities in the random-cluster model we follow \cite{GH96}. For any $ \Delta\subset\L $, $s\in E $, $\bzeta\in\bO_{\L,\bo} $ and $ (\zeta,T) $, with $ T\subset T_\zeta $ we define
 $$
 N_{\Delta,s}(\bzeta)=\#\{\zeta^{\ssup{s}}\cap\Delta\}.
 $$ 
 Then 
 $$
 \begin{aligned}
 N_{\Delta\leftrightarrow\L^{\rm c}}(\zeta,T)=\#\{ &x\in\zeta\cap\Delta\colon x\mbox{ belongs to a cluster connected to } \L^{\rm c} \\ & \mbox{ in } T\cap\Del_3(\zeta)\}
 \end{aligned}
 $$ is the number of points in $ \zeta\cap\Delta $ connected to any point in $ \L^{\rm c} $ in the random hypergraph $T\cap\Del_3(\zeta) $.
 Because of the tile-drawing mechanism, $ \{\Delta\leftrightarrow\L^{\rm c}\} =\{N_{\Delta\leftrightarrow\L^{\rm c}}>0\} $ is also the event that there exists a point in $ \zeta\cap\Delta $ connected to infinity in $ T\cap\Del_3(\zeta) $.
 
 The next Proposition is the key argument why percolation for the random cluster measures leads to a break of symmetry in the mark distribution, and as such it is standard and well known (\cite{GH96}.
 \begin{proposition}\label{Prop-sym}
 For any measurable $\Delta\subset\L $,
 $$
 \int\,(qN_{\Delta,1}-N_\Delta)\,\d\gamma_{\L,\bo} =(q-1)\int\,N_{\Delta\leftrightarrow\L^{\rm c}}\,\d C_{\L,\om}.
 $$
 \end{proposition}
 \begin{proofsect}{Proof}
 This is proved in \cite[Lemma~2.17]{E14} following ideas in  \cite{GH96}.
 \qed
 \end{proofsect}
 
 For purely Delaunay edge potentials an analogous Delaunay random-cluster model exists for edge percolation. Then $ T_{\R^2} $ is replaced by $ E_{\R^2}=\{e=\{x,y\}\subset\R^2\colon x\not= y\} $ and the tile drawing mechanism by an edge drawing mechanism (for details see \cite{BBD04} or \cite{GH96}).

\section{Tile-Percolation}\label{tileperc}
  We establish the existence of tile percolation for the Delaunay random-cluster measure $ C_{\L,\om} $ when $ z$ and the parameter $ \beta $ are sufficiently large. Note that for any $ \Delta\subset\R^2 $ we write
  $$
  N_{\Delta\leftrightarrow\infty}(\zeta,T)=\#\{x\in\zeta\cap\Delta\colon x \mbox{ belongs to an }  \infty-\mbox{cluster of } (\zeta,T\cap T_\zeta)\}.
  $$  
  The key step in our results is the following percolation result.
 \begin{proposition}\label{Prop-per}
  Suppose all the assumptions hold and that $ z $ and $ \beta $ are sufficiently large. Suppose that $ \L $ is a finite union of cells $ \Delta_{k,l} $ defined in \eqref{cell} in Appendix A. Then there exists $\eps>0 $ such that
  $$
  \int\, N_{\Delta\leftrightarrow\infty}\,\d C_{\L,\om}\ge \eps
  $$ for any cell $\Delta=\Delta_{k,l} $, any finite union $ \L $ of cells and any admissible pseudo-periodic boundary condition $ \om \in\O_\L^*$.  
  \end{proposition}

\noindent\textbf{Proof of Proposition~\ref{Prop-per}}: We split the proof in several steps and Lemmata below.  
Our strategy to establish percolation in the Delaunay  random-cluster model is to compare it to site percolation on $ \Z^2 $.  We adapt here the strategy in \cite{GH96} to our cases in the following steps. First we employ a coarse-graining strategy to relate each site $(k,l)\in\Z^2 $ to a cell which is a union of parallelotopes. In order to establish site percolation we need to define when cells are good (open) and when two neighbouring cells are linked once they are open. This link establishes then an open connection in our Delaunay hypergraph structures $ \Del_m(\zeta), m=2,3 $. 

\noindent \textbf{Step 1: Coarse graining.}

\noindent  Let $ \L=\L_n\subset\R^2 $ be the parallelotope  given as the  finite union of cells \eqref{cell} with side length $\ell $, i.e.,
 $$
 \L_n=\bigcup_{(k,l)\in \{-n,\ldots,n\}^2}\Delta_{k,l}\quad \mbox{ and } \Delta_{k,l}=\bigcup_{i,j=0}^8\Delta_{k,l}^{i,j},
 $$ 
  where $ \Delta_{k,l}^{i,j} $ are parallelotopes with side length $ \ell/9 $. The central band of two neighbouring cells $ \Delta_{k,l} $ and $ \Delta_{k+1,l} $ is defined by
$$
  CB_{k:k+1,l}=\big(\bigcup_{i=0}^4\Delta_{k,l}^{4+i,4}\big)\cup\big(\bigcup_{i=0}^4\Delta_{k+1,l}^{i,4}\big).
  $$
  
  \begin{figure}[h!]
\caption{The central band $ \boldsymbol{CB_{k:k+1,l}} $ is the dark shaded strip (union of 10 smaller cells) connecting two neighbouring cells. Some hyperedges are drawn which meet the central band.}\label{fig0}
\includegraphics[scale=0.92]{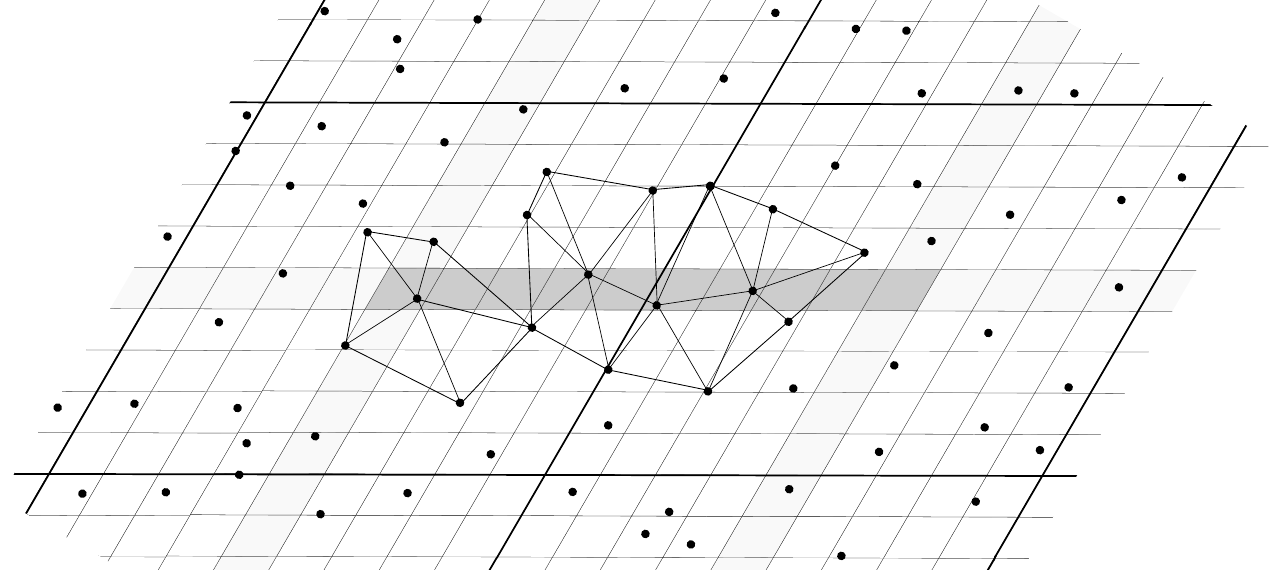}
\end{figure}
 
  \medskip
  
  \noindent We now define two properties our Delaunay random-cluster measures need to satisfy in order to establish percolation in the Delaunay random-cluster models and therefore non-uniqueness of Gibbs measures for certain parameters.  Given $ \zeta\in\O_{\L,\om} $, let $ H_{k:k+1,l}(\zeta) $ be the subset of all tiles of $ \Del_3(\zeta) $ whose circumcircles  have a nonempty intersection with $ CB_{k:k+1,l} $. Define the event $F_{k,l} $ that all small cells $ \Delta_{k,l}^{i,j} $ in $ \Delta_{k,l} $ contain at least one point of $ \zeta $, i.e., $$ F_{k,l}=\bigcap_{i,j=0}^8\{\#\{\zeta\cap\Delta_{k,l}^{i,j}\}\ge 1\}.$$

We say the Delaunay random-cluster measure $ C_{\L,\om} $ is \textbf{coarse grain ready} (\textbf{CGR}) if there are $ 0<\ell_*<\ell^*<\infty $ such that for all scales $ \ell\in(\ell_*,\ell^*) $ there is a continuous, strictly increasing function $g $ of the parameter $ \beta $ such that for all neighbouring cells $ \Delta_{k,l} $ and $ \Delta_{k+1,l} $ and for all configurations $ \zeta\in F_{k,l}\cap F_{k+1,l} $,
$$
 \phi(\alpha(\tau))\ge g(\beta)\qquad \mbox{ for all } \tau\in H_{k:k+1,l}(\zeta).
 $$

 To define the second property we need to consider the  distribution of the points  given by the marginal distribution
 $$
 M_{\L,\om}=C_{\L,\om}(\cdot\times\Tscr)
 $$ of the Delaunay random-cluster measure on $ \O$. Note that  \eqref{DRCmeasure} can be written as
$$
 C_{\L,\om}(\d\zeta,\d T)=M_{\L,\om}(\d\zeta)\mu^{\ssup{q}}_{\zeta,\L}(\d T),\qquad \mu^{\ssup{q}}_{\zeta,\L}(\d T)=\frac{q^{K(\zeta,T)}\mu_{\zeta,\L}(\d T)}{\int\,q^{K(\zeta,T)}\mu_{\zeta,\L}(\d T)}.
$$
  We define $ h_\L $ to be the Radon-Nikodym density of $ M_{\L,\om} $ with respect to $ P_{\L,\om}^z $, i.e., for $ \zeta\in\O_{\L,\om} $,
  $$
  h_\L(\zeta):=\Zscr_{\L}(\om)^{-1}\int\,q^{K(\zeta,T)}\mu_{\zeta,\L}(\d T).
  $$
  We say the Delaunay random-cluster measure $ C_{\L,\om} $ has \textbf{bounded Papangelou conditional intensity} (\textbf{BPI}) if there exists $ \delta=\delta(\psi,\phi)>0 $ such that for $ M_{\L,\om}$-almost all $ \zeta\in\O_{\L,\om} $ and $ x_0\in\L\setminus\zeta $,
  \begin{equation}
  \frac{h_\L(\zeta\cup\{x_0\})}{h_\L(\zeta)}\ge q^{-\delta}.
  \end{equation}

 \medskip 
  
\noindent Suppose our Delaunay random-cluster model satisfies both the (\textbf{CGR})- and the (\textbf{BPI})-property. An important component of our coarse graining method is to estimate  the conditional probability that at least one point lies inside some $ \Delta\subset \L $. In the following we denote by $ M_{\L,\Delta,\om^\prime} $ the conditional distribution (relative to $ M_{\L,\om} $) in $ \Delta $ given that the configuration in $ \L\setminus \Delta $ is equal to $ \om^\prime $ with $ \om^\prime\setminus\L=\om $. The details of the construction of the regular conditional probability distribution can be found in \cite{E14}.
  
  \medskip
  
\noindent Fix $ \eps=\frac{1}{2}(1-p_{\rm c}^{\ssup{\rm site}}(\Z^2)) $. 

\begin{lemma}\label{lem1}
Let $ \L $ as above and let $ \ell>18\delta_0 $. There are constants $ c=c(\delta_0) $ and $ z_0=c(\delta_0)q^\delta $ such that for $ z>z_0 $ and  for all admissible pseudo-periodic boundary conditions $ \om \in\O_\L^*$,
$$
M_{\L,\nabla,\om^\prime}(\#\{\zeta\cap\nabla\}\ge 1)>1-\frac{\eps}{81}
$$
for all cells $ \nabla=\Delta_{k,l}^{i,j}, (k,l)\in\{-n,\ldots,n\}^2, i,j=0,\ldots, 8 $, and for any admissible boundary condition $ \om^\prime\in\O_\nabla^* $ with $ \om^\prime\setminus\L=\om $.
\end{lemma}  
  
\begin{proofsect}{Proof}
The statement follows immediately  from (\textbf{BPI}) and \cite{E14, BBD04} using that the particle distribution is governed by a  hard-core-potential $\psi $ with range $ \delta_0 $. Let $ \nabla_0 $ be the cell $ \nabla $ reduced by a boundary layer of thickness $ \delta_0 $. In particular there is $ z_0$ such that 
\begin{equation}\label{lowerbound}
M_{\L,\nabla,\om^\prime}(\#\{\zeta\cap\nabla\}=0)\le\frac{q^\delta}{z|\nabla_0|}<\frac{\eps}{81}
\end{equation}
for all $ z>z_0 $. For the proof  the coarse graining length scale $ \ell $ must be chosen such that $ \ell>18\delta_0 $.
\qed
\end{proofsect}

  \medskip

  \noindent \textbf{Step 2: Random-cluster measure $\widetilde{C}_{\L,\om}$.}
 We find a measure $ \widetilde{C}_{\L,\om} $ which is stochastically smaller than $ C_{\L,\om} $. Then using coarse graining  and comparison to site percolation on $ \Z^2 $ we establish percolation for $ \widetilde{C}_{\L,\om} $. Percolation for $ \widetilde{C}_{\L,\om} $ then implies percolation for the original $ C_{\L,\om} $. We base the definition of the measure $ \widetilde{C}_{\L,\om} $ on a coarse graining method originally introduced in \cite{Haggstrom00} and then used in \cite{BBD04}.

 \noindent For given $ \zeta\in\O_{\L,\om} $ let 
 $$
 \Del_3^*(\zeta)=\{\tau\in\Del_3(\zeta)\colon \phi(\alpha(\tau))\ge g(\beta)\}
 $$
 and let $ \widetilde{\mu}_{\zeta} $ be the distribution of the random set $\{\tau\in T_\zeta\colon \upsilon(\tau)=1\} $  with
 \begin{equation}\label{siteperc}
 \P(\upsilon(\tau)=1)=\widetilde{p}=\frac{1-\exp\{-g(\beta)\}}{1+(q-1)\exp\{-g(\beta)\}}\1_{\Del_3^*(\zeta)}(\tau).
 \end{equation}
 Note that $ \widetilde{\mu}_\zeta $ depends on $ \L $ only via the configuration $ \zeta\in\O_{\L,\om}$. 
 Note also that $ \widetilde{p} $  is increasing in $\beta$, although to reduce excessive notation, we don't explicitly write this. 
 It is easy to show that $ \mu_{\zeta,\L}\succcurlyeq\widetilde{\mu}_{\zeta} $ (see \cite[Lemma~3.11]{E14}).
 Hence, $ C_{\L,\om}\succcurlyeq\widetilde{C}_{\L,\om} $. As site percolation implies hyperedge percolation we consider site percolation given by \eqref{siteperc}, that is, we open vertices in $ \Del_1(\zeta) $  independently of each other with probability $ \widetilde{p}$. Formally this is defined as follows.  We let $ \widetilde{\lambda}_{\zeta,\L} $ be the distribution of the random mark vector $ \widetilde{\sigma}\in E^\zeta $ where $ (\widetilde{\sigma}_x)_{x\in\zeta} $ are Bernoulli random variables satisfying
 \begin{equation}\label{ptilde}
 \begin{aligned}
 \P(\widetilde{\sigma}_x=1)&=\widetilde{p}\1_{\Del_1^*(\zeta)}(x)\\
 \P(\widetilde{\sigma}_x\not= 1)&=1-\widetilde{p}\1_{\Del_1^*(\zeta)}(x),
 \end{aligned}
 \end{equation}
 where $ \widetilde{p} $ is given in \eqref{siteperc} and $ \Del_1^*(\zeta) $ is the set of points that build the hyperedges of $ \Del_3^*(\zeta) $.
Then the site percolation process is defined by the measure
$$
 \widetilde{C}^{\ssup{\rm site}}_{\L,\om}(\d\bzeta)=M_{\L,\om}(\d\zeta)\widetilde{\lambda}_{\zeta,\L}(\d\bzeta).
 $$ 
 
 \noindent\textbf{Step 3: Site percolation.}
 
Given $ (k,l)\in\{-n,\ldots,n\}^2 $ we need to estimate the probability of the event that all cells $ \Delta_{k,l}^{i,j}, i,j=0,\ldots,8 $, contain at least one point. Pick some $ \nabla=\Delta_{k,l}^{i,j} $ and choose according to Lemma~\ref{lem1} $ z $ so large that for all admissible configurations $ \om^\prime\in\O_\nabla^* $ with $ \om^\prime\setminus\L=\om $, we obtain
$$
M_{\L,\nabla,\om^\prime}(\#\{\zeta\cap\nabla\}=0)\le \frac{\eps}{81}.
$$ Hence, for $ \om^\prime\in\O^*_{\Delta_{k,l}} $ we obtain
$$
M_{\L,\Delta_{k,l},\om^\prime}(\#\{\zeta\cap\Delta_{k,l}\}=0)\le\frac{\eps}{81} \mbox{ and }M_{\L,\Delta_{k,l},\om^\prime}(F_{k,l})>1-\eps.
$$
 
We now establish percolation for the random-cluster measure $ \widetilde{C}_{\L,\om} $.
 
\begin{lemma}\label{LemPercolation}
There exists $ \beta_0>0  $ such that for $ \beta>\beta_0 $   and $ z>z_0 $  where $ z_0 $ is given in Lemma~\ref{lem1}   there is a  $c>0 $ such that 
$$
\widetilde{C}^{\ssup{\rm site}}_{\L,\om}(\{\Delta\longleftrightarrow\L^{\rm c}\})\ge c>0
$$ for any $ \Delta=\Delta_{k,l}\subset\L, (k,l)\in\{-n,\ldots,n\}^2$, and any pseudo-periodic admissible boundary condition $ \om \in\O^*_{\L}$.
\end{lemma}

\begin{proofsect}{Proof}
Pick $ (k,l)\in\{-n,\ldots,n\}^2 $ and consider the event $ C_{k,l} $ that each cell $ \Delta_{k,l}^{i,j} $ has at least one point and all points in $ \Delta_{k,l}\cap\Del_1^*(\zeta) $ are carrying mark $1$,
$$
C_{k,l}=\{\bzeta\in\bO_{\L,\bo}\colon \zeta\in F_{k,l} \mbox{ and } \sigma_{\bzeta}(x)=1 \mbox{ for all } x\in\Delta_{k,l}\cap\Del_1^*(\zeta\}.
$$
A cell $ \Delta_{k,l} $ is declared to be 'good' or 'open' if $ C_{k,l} $ occurs.
Each vertex $ x\in\Del_1^*(\zeta) $ is open with probability $ \widetilde{p} $ (see \eqref{ptilde}). The hard-core background potential $ \psi $ provides an upper bound on the number of points in $ \Delta_{k,l} $ for any configuration $ \zeta $. Thus with
$$
\#\{\Delta_{k,l}\cap\Del_1^*(\zeta)\}\le\frac{(\ell+2\delta_0)^2}{\pi\delta_0^2}=:M
$$
we obtain
$$
\widetilde{C}_{\L,\om}^{\ssup{\rm site}}(C_{k,l})\ge \int\,M_{\L,\Delta_{k,l},\om^\prime}(\d\zeta)\1_{F_{k,l}}(\zeta)\widetilde{p}^M\ge (1-\eps)\widetilde{p}^M.
$$
Using (\textbf{CGR}) we can now find $ \beta_0 $ such that for all $ \beta\ge \beta_0 $ we have
$$
\widetilde{p}\ge (1-\eps)^{\frac{1}{M}}, 
$$ and thus
$$
\widetilde{C}_{\L,\om}^{\ssup{\rm site}}(C_{k,l})\ge (1-\eps)^2>1-2\eps=p_{\rm c}^{\ssup{\rm site}}(\Z^2).
$$ 
The key point is that this property holds for any pair $ (k,l) $ independently from any other pair of box labels. Standard results for site percolation on $ \Z^2 $ now ensure that the probability of a chain of 'good' cells $ \Delta_{m,r} $ from $ \Delta_{k,l}$ to $ \L^{\rm c} $ is strictly positive. To see this note that each site $(k,l)\in\Z^2 $ is open once $ C_{k,l} $ happens. Therefore there is strictly positive probability connecting some cell to the boundary of $ \L $ via open sites ('good' cells). It remains to check whether this chain also ensures  an open connection (path)  from $ \Delta_{k,l} $ to $ \L^{\rm c} $ within the Delaunay graph $ \Del_3(\zeta)$. Our previous definition of the central band is designed in such a way that two neighbouring cells which are both open ('good') guarantee an open path using hyperedges which have an intersection with the central band.   Suppose $ \bzeta\in C_{k,l}\cap C_{k+1,l} $, that is, the neighbouring cells are 'good'. Since $ C_{\L,\om} $ satisfies (\textbf{CGR}) it follows that $ \phi(\alpha(\tau))\ge g(\beta) >0 $ for all $ \tau\in H_{k:k+1,l}(\zeta) $. Thus $ H_{k:k+1,l}(\zeta)\subset\Del_3^*(\zeta) $ and henceforth $ \sigma_{\bzeta}(x)=1 $ for all $ x\in H_{k:k+1,l}(\zeta) $. In other words, for all $ (k,l)\in\{-n,\ldots,n\}^2 $ the restriction of $ \Del_3^*(\zeta) $ to $ H_{k:k+1,l}(\zeta) $ equals the restriction of $ \Del_3(\zeta) $ to $ H_{k:k+1,l}(\zeta) $. Hence we can connect $ \Delta_{k,l}^{4,4} $ to $ \Delta_{k+1,l}^{4,4} $ in the graph $ \Del_3(\zeta) $ inside $ \Delta_{k,l}\cup\Delta_{k+1,l} $ (via central band hyperedges) connecting points of mark $1$. Hence, there is $c>0 $ with $ \widetilde{C}_{\L,\om}(\{\Delta_{k,l}\longleftrightarrow \L^{\rm c}\})>c>0 $.\qed
\end{proofsect}

\section{Proofs}\label{proofs} 
This section is devoted to the proof of our results. We first establish the existence of a Gibbs measure in our models.
\subsection{Existence of Gibbs measures}
To show the existence of Gibbs measures (Proposition~\ref{Propexist}) for all our Delaunay Potts models we will show that they satisfy  adapted versions of the conditions necessary for the existence outlined in \cite{DDG12}. For Delaunay hypergraph  structures the authors of \cite{DDG12} have given three general conditions hyperedge potentials need to satisfy to ensure existence of Gibbs measures.  The potentials $ \phi $ and $ \psi $ of all our models depend solely on the Delaunay hyperedges $ \Del_m(\bzeta) $, $m=2 $ or $ m=3 $,  of a marked configuration $ \bzeta $. Every marked hyperedge $ \beeta\in\Del_2(\bzeta) $ has the so-called finite horizon $\overline{B}(\eta,\zeta) $ where $ B(\eta,\zeta)$ is the open ball with $ \partial B(\eta,\zeta)\cap\zeta=\eta $ that contains no points of $ \zeta$. The finite-horizon property of a general hyperedge potential $ \varphi\colon\bO_f\times\bO\to\R  $ says that for each pair $ (\beeta,\bzeta) $ with $ \beeta\in\Del_m(\bzeta) $ there  exists some $ \Delta\Subset\R^2 $ such that for the pair $ (\beeta,\widetilde{\bzeta}) $ with $ \beeta\in\Del_m(\widetilde{\bzeta}) $ we have that $ \varphi(\beeta,\bzeta)=\varphi(\beeta,\widetilde{\bzeta}) $ when $ \widetilde{\bzeta}=\bzeta $ on $ \Delta\equiv\overline{B}(\eta,\zeta) $. The hyperedge potential $ \varphi $ in all our models is the sum of two potentials, the background potential $ \psi $ and the type potential $ \phi $. As all our potentials depend solely on the single hyperedge they share this finite-horizon property.  Thus they satisfy the range condition (R) in \cite{DDG12}, see \cite[Proposition~4.1 \& 4.3]{DDG12}. The second requirement is the stability condition. A hyperedge potential  $ \varphi $ is called stable if there is a lower bound for the Hamiltonian for any $ \L\Subset\R^2 $. This stability condition (S) is satisfied for all our potentials $\psi $ and $ \phi $ as they are all positive. The third condition to be checked  is a partial complementary upper bound for the Hamiltonian in any $ \L\Subset \R^2 $. This is a bit more involved and we shall first define appropriate configurations, the so-called pseudo-periodic marked configurations. We consider the partition of $ \R^2 $ as given in Appendix~\ref{pseudoperiodic}. We let $ B(0,r) $ be an open ball of radius $ r\le \rho_0\ell $ where we choose $ \rho_0>0 $ sufficiently small such that  $ B(0,r)\subset\Delta_{0,0} $. Note that 
$$ B^r:=\{\zeta\in\O_{\Delta_{0,0}}\colon \zeta=\{x\}\mbox{ for some } x\in B(0,r)\} $$ is a measurable set of $\O_{\Delta_{0,0}}\setminus\{\emptyset\} $. Then
$$
\Gamma^r=\{\om\in\O\colon\theta_{Mz}(\om_{\Delta_{k,l}})\in B^r\mbox{ for all } (k,l)\in\Z^2\} 
$$ is the set of pseudo-periodic configurations \eqref{ppc}. Note that these configurations are not marked. The reason is that when a point is shifted its mark remains unchanged. Thus we define the set of pseudo-periodic marked configurations as 
$$
\boldsymbol{\Gamma}^r=\{\bo=(\om^{\ssup{1}},\ldots,\om^{\ssup{q}})\colon \om^{\ssup{i}}\in\Gamma^r\mbox{ for all }i\in E\}.
$$

The required control of the Hamiltonian from above will be achieved by the following properties. As all our hyperedge potentials depend only on the single hyperedge the so-called uniform confinement (see \cite{DDG12}) is trivially satisfied. In addition we need the uniform summability for any hyperedge potential, that is,
$$
c_r:=\sup_{\bzeta\in\boldsymbol{\Gamma}^r}\sum_{\beeta\in\Del_m(\bzeta)\colon \eta\cap\Delta\not=\emptyset}\frac{\varphi(\beeta,\bzeta)}{\#\widehat{\eta}}<\infty,
$$
where $ \widehat{\eta}=\{(k,l)\in\Z^2\colon \eta\cap \Delta_{k,l}\not=\emptyset\} $. Note that the length $ \ell(\eta) $ of any $ \eta\in\Del_2(\zeta)\cap \Delta $ when $ \zeta $ is any pseudo-periodic configuration  satisfies  $ \ell(1-2\rho_0)\le \ell(\eta)\le \ell(1+2\rho_0) $.  As all our models have the background hard-core potential $ \psi $ with range $ \delta_0 $ we need to choose  the scale $ \ell $ for the partition such that $ \ell(1-2\rho_0) > \delta_0 $. 
We compute the constant $ c_r $ for each model. Our chosen partition  implies that maximal six Delaunay hyperedges exists. Furthermore, $ \widehat{\eta}=2 $ when $ \eta\in\Del_2(\zeta) $  and   $ \widehat{\eta}=3  $ when $ \eta\in\Del_3(\zeta) $.

\noindent Triangle model I: As $ \phi(\alpha(\tau))\le \beta $ for any $ \tau\in\Del_3(\zeta) $ we obtain
$$
c_r=2\beta,
$$
and thus uniform summability holds as long as $ \beta<\infty $. The so-called strong rigidity is given once $ z|B(0,r)|>\ex^{c_r} $, that is, once
$$
z>\frac{\ex^{2\beta}}{\pi\rho_0^2(18\delta_0)^2}=:z_0(\delta_0,\beta).
$$

\noindent Triangle model II: As $ \phi(\alpha(\tau))\le \log(1+\beta(\pi/3)^3) $ for any $ \tau\in\Del_3(\zeta) $ we obtain
$$
c_r=2\log(1+\beta(\pi/3)^3)$$
and thus uniform summability. We obtain strong rigidity once
$$
z>\frac{1}{\pi\rho_0^2(18\delta_0)^2}\big(1+(\beta\pi^3)/27\big)^{2}=:z_0(\delta_0,\beta).
$$
\noindent Edge model: 
As $ \phi(\ell(\eta))\le \log(1+\beta(\delta_0/\ell(1-2\rho_0))^3) $ we obtain
$$
c_r=3\log\Big(1+\frac{\beta}{(18(1-2\rho_0))^3}\Big)<\infty,
$$ and thus uniform summability. We obtain strong rigidity once 
$$
z>\frac{1}{\pi\rho_0^2(18\delta_0)^2} \Big(1+\frac{\beta}{(18(1-2\rho_0))^3}\Big)^3    =:z_0(\delta_0,\beta).
$$
Using \cite[Theorem~3.3]{DDG12} and \cite[Corollary~3.4]{DDG12} we obtain all the statements in Proposition~\ref{Propexist} by noting that for any given $ z>0 $ we can pick $ \rho_0 $ and $ \ell $ sufficiently large such that the above conditions on $ z $ are satisfied. \qed

\subsection{Percolation of Delaunay random cluster models}
We show percolation in all Delaunay random-cluster models for sufficiently large parameter $ z $ and $ \beta $. As outlined in Section~\ref{tileperc} we need to establish the (\textbf{CGR})-property and the (\textbf{BPI})-property.

\begin{lemma}[\textbf{CGR - Triangle Model I}]\label{LemTMI}
For all $ \zeta\in F_{k,l}\cap F_{k+1,l} $ and $ \tau\in H_{k:k+1,l}(\zeta) $, $ (k,l)\in\{-n,\ldots,n\}^2 $,
$$
\phi(\alpha(\tau))\ge \beta,
$$ where the coarse graining scale $ \ell $ satisfies $ \ell\in (18\delta_0,m\delta_0) $ for some $ m>18 $  and $ \alpha_0\in (0,\frac{9}{m\sqrt{7}}) $.
\end{lemma}
\begin{proofsect}{Proof}
Suppose that $ 2\delta_0<\ell< m\delta_0$ and $ \alpha_0\in (0,\frac{9}{m\sqrt{7}}) $. Since all of the little cells $ \Delta_{k,l}^{i,j}, i,j=0,\ldots, 8 $ with side length $ \ell/9 $, contain at least one point, we have that any open ball of radius at least $ \frac{1}{18}\sqrt{7}\ell $  and centre $ x_0\in CB_{k:k+1,l} $ has a non-empty  intersection with $ \zeta $. Therefore, for each $ \tau\in H_{k:k+1,l}(\zeta) $  the circumscribing circle $ B(\tau,\zeta) $ has radius less than $ \frac{1}{18}\sqrt{7}\ell $. Let $ \tau=\{a,b,c\} $ be such a triangle. Without loss of generality, let $ \alpha(\tau) $ be the angle $ \widehat{acb} $ and $ l$ be the arc length of the arc on $ \partial B(\tau,\zeta) $ between $a $ and $b$. Let $ x$ be the centre of $ B(\tau,\zeta) $. It follows that $ \widehat{axb}=2\alpha(\tau) $ and $ l=2r\alpha(\tau) $ where $ r $ is the radius of $ B(\tau,\zeta) $. By the hard-core condition $ l>|a-b|\ge \delta_0 $. Combing these facts with $ r<\frac{1}{18}\sqrt{7}\ell $ gives
$$
\alpha(\tau)\ge \frac{9\delta_0}{\sqrt{7}\ell}\ge \frac{9}{\sqrt{7}m}>\alpha_0,
$$ and thus $\phi(\alpha(\tau))\ge \beta $ for all  $\tau\in H_{k:k+1,l}(\zeta)$. We obtain (\textbf{CGR}) by choosing $ g(\beta)=\beta $.
\qed
\end{proofsect}

\begin{lemma}[\textbf{CGR - Triangle Model II}]\label{CGRII}
For all $ \zeta\in F_{k,l}\cap F_{k+1,l} $ and $ \tau\in H_{k:k+1,l}(\zeta) $, $ (k,l)\in\{-n,\ldots,n\}^2 $,
$$
\phi(\alpha(\tau))\ge \log\big(1+\beta\alpha_0^3\big),
$$ where the coarse graining scale $ \ell $ satisfies $ \ell\in(18\delta_0,m\delta_0) $ for some $ m>18$.
\end{lemma}
\begin{proofsect}{Proof}
With $ \alpha(\tau)\ge  \alpha_0 $ from Lemma~\ref{LemTMI} we obtain the statement. Note that the lower bound is strictly increasing in the parameter $ \beta $. Thus we conclude with $ g(\beta)=\log\big(1+\beta\alpha_0^3\big)$.\qed
\end{proofsect}

\begin{lemma}[\textbf{CGR - Edge Model}]\label{CGRIII}
For all $ \zeta\in F_{k,l}\cap F_{k+1,l} $ and all edges $ \eta\in\Del_2(\zeta) $ such that $ \eta\subset\tau $ for some $ \tau\in H_{k:k+1,l}(\zeta) $, $ (k,l)\in\{-n,\ldots,n\}^2 $,
$$
\phi(\ell(\eta))\ge \log\big(1+\beta/(9\sqrt{7}m)^3\big),
$$ where the coarse graining scale $ \ell $ satisfies $ \ell\in(18\delta_0,m\delta_0) $ for some $ m>18 $.
\end{lemma}
\begin{proofsect}{Proof}
From the proof of Lemma~\ref{LemTMI} we know that the radius of the circumcircle of any triangle $ \tau\in H_{k:k+1,l}(\zeta) $ is less than $ \frac{1}{18}\sqrt{7}\ell $ and thus any Delaunay edge in $ \Del_2(\zeta) $ is no longer than $ \frac{1}{9}\sqrt{7}\ell $. We thus conclude with the lower bound noting that this bound is strictly increasing in $ \beta $ by choosing $ g(\beta)=\log\big(1+\beta/(9\sqrt{7}m)^3\big)$.\qed
\end{proofsect}

\medskip

\begin{figure}[h!]
\caption{The hyperedge sets $ T_{x_0,\zeta}^{\ssup{\rm ext}}, T^{\ssup{+}}_{x_0,\zeta} $, and $ T^{\ssup{-}}_{x_0,\zeta} $ from the left to right}\label{fig}
\includegraphics[scale=0.58]{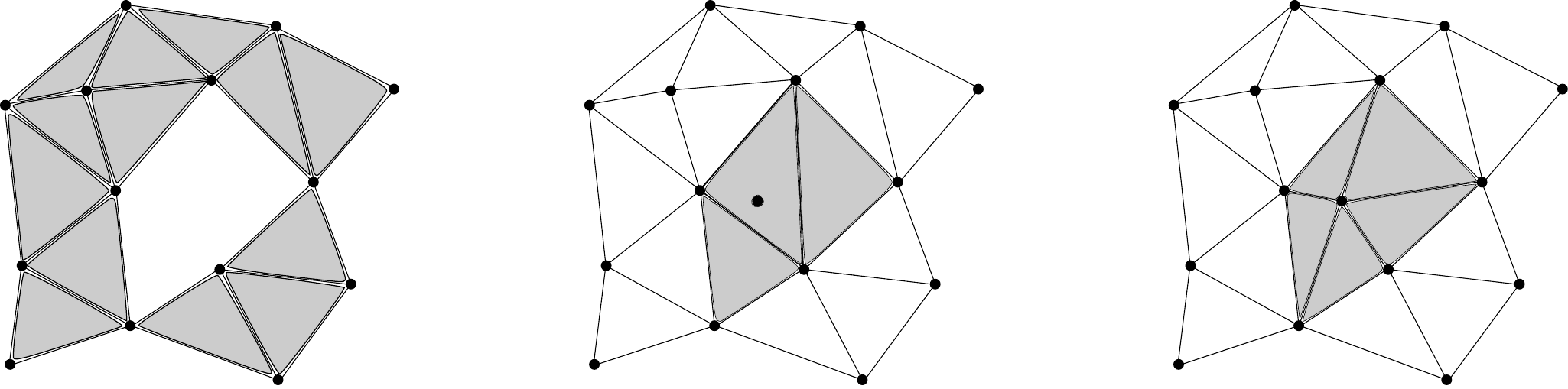}
\end{figure}

In order to show (\textbf{BPI}) for all our models, we need to investigate the geometry of the Delaunay triangulation $ \Del_3(\zeta) $, and in particular what happens to it when we augment $ \zeta $ with a new point $ x_0\notin\zeta $. Some hyperedges may be destroyed, some are created, and some remain. This process is well described in \cite{Lis94}. We give a brief account here for the convenience of the reader. We insert the point $ x_0 $ into one of the triangles $ \tau $ in $ \Del_3(\zeta) $. We then create three new edges that join $ x_0 $ to each of the three vertices of $ \tau $. This creates three new triangles, and destroys one. We now need to verify that the new triangles each satisfy the Delaunay condition \eqref{triangulation}, that is, that their circumscribing balls contain no points of $ \zeta$. If this condition is satisfied the new triangle remains, if it is not satisfied, then there is a point $ x_1\in\zeta $  inside the circumscribing ball. We remove the edge not connected to $ x_0 $, and replace it by an edge connecting $x_0 $ and $ x_1 $. This results in the creation of two new triangles. Each of these triangles must be checked as above and the process continues. Once all triangles satisfy the Delaunay condition, we arrive at the Delaunay triangulation $ \Del_3(\zeta\cup\{x_0\}) $. Let
\begin{equation}\label{tilesets}
\begin{aligned}
T_{x_0,\zeta}^{\ssup{\rm ext}}&=\Del_3(\zeta)\cap\Del_3(\zeta\cup\{x_0\}),\\
T^{\ssup{+}}_{x_0,\zeta}&=\Del_3(\zeta\cup\{x_0\})\setminus\Del_3(\zeta)=\Del_3(\zeta\cup\{x_0\})\setminus T_{x_0,\zeta}^{\ssup{\rm ext}},\\
T^{\ssup{-}}_{x_0,\zeta}&=\Del_3(\zeta)\setminus\Del_3(\zeta\cup\{x_0\})=\Del_3(\zeta)\setminus T_{x_0,\zeta}^{\ssup{\rm ext}},
\end{aligned}
\end{equation}
be the set of exterior, created, and destroyed triangles respectively, see figure~\ref{fig}. 

Note that any new triangle must contain $ x_0 $, i.e.,
$$
T^{\ssup{+}}_{x_0,\zeta}=\{\tau\in\Del_3(\zeta\cup\{x_0\})\colon \tau\cap x_0=x_0\}.$$
We let $ \mu^{\ssup{-}}_{x_0,\zeta},\mu^{\ssup{+}}_{x_0,\zeta} $, and $ \mu^{\ssup{\rm ext}}_{x_0,\zeta} $ be the tile (triangle) drawing mechanisms on $ T_{x_0,\zeta}^{\ssup{\rm ext}}, T^{\ssup{+}}_{x_0,\zeta} $, and $ T^{\ssup{-}}_{x_0,\zeta}$,  respectively, which are derived from the tile drawing measure $ \mu_{\zeta,\L} $ in Section~\ref{DelCluster}. We need to study the change of $ K(\zeta,T), T\subset T_\zeta\subset\Tscr $ when adding the point $ x_0\notin\zeta $. Adding a point $x_0 $ to $ \zeta $ without considering the change to $ T $ will always increase the number of connected components by one. On the other hand, the augmentation of a single triangle $ \tau $ to $T$ can result in the connection of a maximum of three different connected components, leaving one. Therefore,
\begin{equation}\label{change}
\begin{aligned}
K(\zeta\cup\{x_0\},T)-K(\zeta,T)&=1,\\
-2\le K(\zeta,T\cup\tau)-K(\zeta,T)&\le 0.
\end{aligned}
\end{equation}

\begin{lemma}[\textbf{BPI - Triangle Model I}]\label{LemBPI}
For $ \L\Subset\R^2 $ and pseudo-periodic admissible boundary conditions  $\om\in\O_{\L}^* $, and $ M_{\L,\om} $-almost all $ \zeta\in\O_{\L,\om} $ and $ x_0\in\zeta\setminus\L $,
$$
\frac{h_\L(\zeta\cup\{x_0\})}{h_\L(\zeta)}\ge q^{-\frac{4\pi}{\alpha_0}}.
$$

\end{lemma}
\begin{proofsect}{Proof}
The tile drawing mechanism $ \mu_{\zeta,\L} $ opens only triangles $ \tau $  with minimal angle $\alpha(\tau)\ge \alpha_0 $. We collect all triangles touching $ \L $ and having minimal angle greater than $ \alpha_0 $ into the set $T^*$. Using this fact and $ \Del_3(\zeta\cup\{x_0\})=T^{\ssup{\rm ext}}_{x_0,\zeta}\cup T^{\ssup{+}}_{x_0,\zeta} $, we get
$$
\begin{aligned}
& \frac{h_\L(\zeta\cup \{x_0\})}{h_\L(\zeta)}=\frac{\int_{T^*}q^{K(\zeta\cup\{x_0\},T)}\mu_{\zeta\cup\{x_0\},\L}(\d T)}{\int_{T^*}q^{K(\zeta,T)}\mu_{\zeta,\L}(\d T)}\\
&= \frac{\int_{T^*\cap T^{\ssup{\rm ext}}_{x_0,\zeta}}\int_{T^*\cap T^{\ssup{+}}_{x_0,\zeta}}q^{K(\zeta\cup\{x_0\},T_1\cup T_2)-K(\zeta,T_1)}\mu^{\ssup{+}}_{x_0,\zeta}(\d T_2)\mu^{\ssup{\rm ext}}(\d T_1)}{\int_{T^*\cap T^{\ssup{\rm ext}}_{x_0,\zeta}}\int_{T^*\cap T^{\ssup{-}}_{x_0,\zeta}}q^{K(\zeta, T_3\cup T_4)-K(\zeta,T_3)}\mu^{\ssup{-}}_{x_0,\zeta}(\d T_4)\mu^{\ssup{\rm ext}}_{x_0,\zeta}(\d T_3)}.
\end{aligned}
$$ Note that the maximal number of triangles in $ T_2 $ is $ \frac{2\pi}{\alpha_0} $. Therefore, by \eqref{change}, we conclude with
$$
\begin{aligned}
K(\zeta\cup\{x_0\},T_1\cup T_2)-K(\zeta,T_1)&\ge -\frac{4\pi}{\alpha_0},\\
K(\zeta,T_3\cup T_4)-K(\zeta,T_3)&\le 0.
\end{aligned}
$$
\qed
\end{proofsect}

\noindent For our triangle model II we shall obtain a lower bound for the difference
$$
K(\zeta\cup\{x_0\},T_1\cup T)-K(\zeta,T_1),
$$ where $ T_1\subset T^{\ssup{\rm ext}}_{x_0,\zeta} $ and $ T\subset T^{\ssup{+}}_{x_0,\zeta} $. Application of \eqref{change} would result in
\begin{equation}\label{diff} 
K(\zeta\cup\{x_0\},T_1\cup T)-K(\zeta,T_1)\ge -2 \#\{T\}.
\end{equation}
For the triangle model I it was easy to bound the number of elements in $T$ via the minimal angle condition. In the triangle model II we only obtain a bound on the expectation of that number.  We obtain the bound on the expectation.

\begin{lemma}\label{expectedtiles}
Let $ \L\Subset\R^2 $ and $\zeta\in\O_{\L,\om} $. Then
$$
\int\,\#\{T\}\,\mu^{\ssup{+}}_{x_0,\zeta}(\d T)\le 2\pi(1+\frac{1}{3}\beta\pi^2)
$$ for $ x_0\in\L\setminus\zeta $.
\end{lemma}
\begin{proofsect}{Proof}
As any $ \tau\in T^{\ssup{+}}_{x_0,\zeta} $ touches $ \L $ we get with the tile drawing probability $ p(\tau) $ in \eqref{tiledrawing},
$$
\int\,\#\{T\}\,\mu^{\ssup{+}}_{x_0,\zeta}(\d T)=\sum_{\tau\in T^{\ssup{+}}_{x_0,\zeta}}\big(1-\ex^{-\phi(\alpha(\tau))}\big).
$$
Write $ T^{\ssup{+}}_{x_0,\zeta} $ as a union of the disjoint sets
$$
\begin{aligned}
H_1&=\{\tau\in T^{\ssup{+}}_{x_0,\zeta}\colon \alpha(\tau)\ge 1\},\\
H_k&=\big\{\tau\in T^{\ssup{+}}_{x_0,\zeta}\colon \frac{1}{k}\le\alpha(\tau)<\frac{1}{k-1}\big\}, k\ge 2.
\end{aligned}
$$ 
Note that $ \#\{H_k\}\le 2\pi k $ for all $ k\in\N $. Thus
$$
\begin{aligned}
\int\,\#\{T\}&\,\mu^{\ssup{+}}_{x_0,\zeta}(\d T)=\sum_{k\in\N}\sum_{\tau\in H_k}\big(1-\ex^{-\phi(\alpha(\tau))}\big)\\& \le 2\pi\big(1+\beta\sum_{k= 2}^\infty\frac{k}{(k-1)^3+1}\big)\le 2\pi(1+\frac{1}{3}\beta\pi^2),
\end{aligned}
$$
where we used that 
$$
\sum_{k=2}^\infty\frac{k}{(k-1)^3}\le 2\sum_{k=1}^\infty\frac{1}{k^2}=\frac{\pi^2}{3}.
$$
\qed
\end{proofsect}

We then get (\textbf{BPI}) using Jensen's inequality. 

\begin{lemma}[\textbf{BPI - Triangle Model II}]
For $ \L\Subset\R^2 $ and pseudo-periodic admissible boundary conditions  $\om\in\O_\L^* $, and $ M_{\L,\om} $-almost all $ \zeta\in\O_{\L,\om} $ and $ x_0\in\zeta\setminus\L $,
$$
\frac{h_\L(\zeta\cup\{x_0\})}{h_\L(\zeta)}\ge q^{-4\pi(1+\frac{1}{3}\beta\pi^2)}.
$$

\end{lemma}
\begin{proofsect}{Proof}
We proceed as in Lemma~\ref{LemBPI}. Jensen's inequality, \eqref{diff}, and Lemma~\ref{expectedtiles} provide the lower bound 
$$
\begin{aligned}
\int\,q^{K(\zeta\cup\{x_0\},T_1\cup T)-K(\zeta,T_1)} & \mu^{\ssup{+}}_{x_0,\zeta}(\d T)\ge \int\,q^{-2\#\{T\}}\mu^{\ssup{+}}_{x_0,\zeta}(\d T)\\&\ge q^{-2\int\#\{T\}\mu^{\ssup{+}}_{x_0,\zeta}(\d T)}\ge q^{-4\pi(1+\frac{1}{3}\beta\pi^2)}.
\end{aligned}
$$
\qed
\end{proofsect}

\noindent For the edge model we define Delaunay  edge sets  $ E^{\ssup{\rm ext}}_{x_0,\zeta}, E^{\ssup{+}}_{x_0,\zeta}, $ and $ E^{\ssup{-}}_{x_0,\zeta} $ in the same way as done for triangles in \eqref{tilesets} (see also \cite{BBD04}).  As outlined in Section~\ref{DelCluster} an edge drawing measure $ \nu_{\zeta,\L} $ is defined in the same way as the tile drawing measure $ \mu_{\zeta,\L} $ in \eqref{tiledrawing}, that is,
$$
p_\L(\eta)=\begin{cases} \big(1-\ex^{\phi(\ell(\eta))}\big)\1_{\Del_2(\zeta)}(\eta) & \mbox{ if } \eta\in E_{\R^2}\setminus E_{\L^{\rm c}},\\
\1_{\Del_2(\zeta)}(\eta) & \mbox{ if } \eta\in E_{\L^{\rm c}},
\end{cases}
$$ where $ \ell(\eta) $ is the length of the Delaunay edge $ \eta$. Then $ \nu^{\ssup{+}}_{x_0,\zeta} $ denotes the corresponding edge drawing measure on $ E^{\ssup{+}}_{x_0,\zeta} $.
We need to bound the expectation of the  edges in $ E^{\ssup{+}}_{x_0,\zeta} $.

\begin{lemma}\label{expectededges}
Let $ \L\Subset\R^2 $ and $\zeta\in\O_{\L,\om} $. Then
$$
\int\,\#\{E\}\,\nu^{\ssup{+}}_{x_0,\zeta}(\d E)\le 4\big(9+\frac{1}{6}\beta\pi^2\big)
$$ for $ x_0\in\L\setminus\zeta $.
\end{lemma}
\begin{proofsect}{Proof}
Due to the hard-core background potential we obtain almost surely with respect to $ M_{\L,\om} $ the following upper bound on the number of points $ y $ in $ \zeta $ with $ |y-x_0|<r $,
$$
|B(x_0,r)\cap\zeta|\le 4\big(\frac{r+\delta_0}{\delta_0}\big)^2,
$$ and thus for $ r=2\delta_0,4\delta_0,\ldots $,
$$
\begin{aligned}
|B(x_0,2\delta_0)\cap \zeta|&\le 36=:b_1,\\
|B(x_0,2n\delta_0)\cap\zeta|&\le 4(2n+1)^2=:b_n, n\ge 1.
\end{aligned}
$$
The number of edges in $ E^{\ssup{+}}_{x_0,\zeta} $ that have length less than $ r $ is bounded by $ |B(x_0,r)\cap\zeta| $, however, as $ r$ increases, this bound grows quadratically. On the other hand, according to the edge potential $ \phi $, edges with large length are less likely to survive the $ p_\L $-thinning process of $ E^{\ssup{+}}_{x_0,\zeta} $ than their counterparts with small length. We tradeoff these facts in the following.
Let
$$
\R^2=\bigcup_{n=0}^\infty A_n\mbox{ with } A_0=B(x_0,\delta_0) \mbox{ and } A_n=B(x_0,2n\delta_0)\setminus A_{n-1}, n\ge 2.
$$ Now let $E$ be the $p_\L $-thinning of $E^{\ssup{+}}_{x_0,\zeta} $ and define
$$
E^n=\{\eta=\{x_0,x\}\in E\colon |x_0-x|<2n\delta_0\}.
$$ Clearly, $ \int\#\{E^1\}\nu^{\ssup{+}}_{x_0,\zeta}(\d E)\le 36 $. Now for $n=2 $ at most $36$ points are in distance less than $2\delta_0 $ from $ x_0 $, whereas the remaining points lie in the annulus $A_2 $, and due to the fact that $p_\L $ is a decreasing function of the length (distance), they have at most probability $ p_\L(2\delta_0) $ of sharing an edge with $ x_0 $ in $E$. Therefore,
$$
\int\,\#\{E^n\}\,\nu^{\ssup{+}}_{x_0,\zeta}(\d E)\le 36+\sum_{k=1}^{n-1} (b_{k+1}-b_k)p_\L(2k\delta_0)\le 4\big(9+8\sum_{k=1}^{n-1}kp_\L(2k\delta_0)\big).
$$
Using
$$
p_\L(k\delta_0)=1-\ex^{-\phi(k\delta_0)}=\frac{\beta}{k^3+\beta},
$$ and the monotone convergence theorem, we finally obtain the bound 
$$
\begin{aligned}
\int\,\#\{E\}\,\nu^{\ssup{+}}_{x_0,\zeta}(\d E)\le4\big(9+\beta \sum_{k=1}^\infty\frac{1}{k^2}\big)=4\big(9+\frac{1}{6}\beta\pi^2\big).
\end{aligned}
$$
\qed
\end{proofsect}
 
 Having this bound we immediately obtain the (\textbf{BPI}) for the  edge model using Jensen's inequality.
 
\begin{lemma}[\textbf{BPI - Edge Model}]\label{BPIIII}
For $ \L\Subset\R^2 $ and pseudo-periodic admissible boundary conditions  $\om\in\O_\L^* $, and $ M_{\L,\om} $-almost all $ \zeta\in\O_{\L,\om} $ and $ x_0\in\zeta\setminus\L $,
$$
\frac{h_\L(\zeta\cup\{x_0\})}{h_\L(\zeta)}\ge q^{-4(9+\frac{1}{6}\beta\pi^2)}.
$$

\end{lemma}

 \noindent We finish the proof of percolation in all Delaunay random-cluster models:\\[1ex]
 \textbf{Triangle Model I:}  We notice that (\textbf{BPI}) holds for $ \delta=\frac{4\pi}{\alpha_0} $ (Lemma~\ref{LemBPI}). We now pick $ z_0=z_0(\delta_0,\alpha_0) $ according to Lemma~\ref{lem1} and  from (\textbf{CGR}) in Lemma~\ref{LemTMI} we obtain $ \beta_0=\beta_{0}(\delta_0) $ such that Lemma~\ref{LemPercolation} holds for $ z>z_0 $ and $ \beta>\beta_0$. Henceforth, we have that Proposition~\ref{Prop-per} holds and with Section~\ref{finish} we conclude with Theorem~\ref{THM-main}(a).
 
 \medskip
 \noindent \textbf{Triangle Model II:} We see that (\textbf{BPI}) holds for $ \delta=4\pi(1+\frac{1}{3}\beta\pi^2) $ and observe that it depends on $ \beta $.  Thus we first pick $ \beta_0 $ such that  $ \widetilde p\ge (1-\eps)^{\frac{1}{M}} $ in the proof of Lemma~\ref{LemPercolation}. We then pick $ z_0=z_0(\delta_0,\beta) $ according to Lemma~\ref{lem1}, i.e., \eqref{lowerbound}, and  from (\textbf{CGR}) in Lemma~\ref{CGRII} we obtain  Lemma~\ref{LemPercolation} such that  Proposition~\ref{Prop-per} is established for all $z>z_0 $,  and finally Theorem~\ref{THM-main}(b).
 
 \medskip
 
 \noindent \textbf{Edge Model:} We see that (\textbf{BPI}) holds with $ \delta=4(9+\frac{1}{6}\beta\pi^2) $ and observe that it depends on $ \beta $. Thus we proceed as before for the triangle model II. Note that we need an adaptation of Proposition~\ref{Prop-per} to the hypergraph structure $ \Del_2(\zeta) $. This follows as outlined in \cite{E14} or \cite{BBD04}. We thus have established Theorem~\ref{THM-mainedge}.

\subsection{Breaking of the symmetry of the mark distribution}\label{finish}
In this section we complete the proof of Theorem~\ref{THM-mainedge} and Theorem~\ref{THM-main} by analysing the Gibbs distributions $ \gamma_{\L,\bo} $ in the limit $ \L\uparrow\R^2 $. We pick an admissible boundary condition $ \om\in\O^*_{\L_n} $ and let $ \bo=(\om\setminus\L_n,\emptyset,\ldots,\emptyset) $ be the admissible monochromatic boundary condition such that $ \bo\in\bO_{\L_n}^* $. We write $ \gamma_n $ for $ \gamma_{\L_n,\bo} $ for ease of notation and we let $ P_n $ be the probability measure on $ \bO $ relative to which the marked configurations in distinct parallelotopes $ \L_n+(2n+1)M(k,l), (k,l)\in\Z^2 $, are independent with identical distribution $ \gamma_n $. As we are dealing with a cell structure for the partition of $ \R^2 $ we confine ourself first to lattice shifts when we employ spatial averaging. Thus,
$$
\overline{P}_n=\frac{1}{2n+1}\sum_{(k,l)\in\{-n,\ldots,n\}^2}P_n\circ\theta_{M(k,l)}^{-1}.
$$
By the periodicity of $ P_n $ the measure $\overline{P}_n $ is $\Z^2$-shift-invariant. The proof in \cite[Chapter 5]{DDG12} shows that $ (\overline{P}_n)_{n\ge 1} $ has a subsequence which converges with respect to the topology of local convergence to some $ \widehat{P} \in\Mcal_1(\bO) $. As outlined in \cite{DDG12} it is difficult to show that $ \widehat{P} $ is concentrated on admissible configurations. As $ \widehat{P} $ is non-degenerate the proof in \cite[Chapter 5]{DDG12} shows that $ P=\widehat{P}(\cdot|\{\emptyset\}^{\rm c}) $ is a Gibbs measure with $ P(\{\emptyset\})=0 $.  In order to obtain an $ \R^2$-shift-invariant Gibbs measure one needs to apply another averaging,
$$
P^{\ssup{1}}=\int_{\Delta_{0,0}}\,P\circ\theta_{Mx}^{-1}\,\d x.
$$
Applying Propositions~\ref{Prop-sym} and \ref{Prop-per}, we see that for $ \Delta=\Delta_{0,0} $,
$$
\begin{aligned}
\int\,(qN_{\Delta,1}-N_{\Delta})\,\d \overline{P}_n& \ge \frac{(q-1)}{2n+1}\sum_{(k,l)\in\{-n,\ldots,n\}^2}\int\,N_{\Delta_{k,l}\leftrightarrow\infty}\,\d C_{\L_n,\om}\\
&\ge (q-1)\eps.
\end{aligned}
$$ Thus
$$
\int\,(qN_{\Delta,1}-N_\Delta)\,\d P^{\ssup{1}}>0,
$$ and using the symmetry of $ P^{\ssup{1}} $ we observe the following break of symmetry in the expected density of particles of type $1 $ and of any other different type, i.e.,
$$
\rho_1(P^{\ssup{1}})>\rho_2(P^{\ssup{1}})=\cdots=\rho_q(P^{\ssup{1}}),
$$ where $ \rho_s(P^{\ssup{1}})=1/|\Delta|\E_{P^{\ssup{1}}}[N_{\Delta,s}], s\in E $. We conclude with our statement as in \cite{GH96} by showing that the matrix 
$$
\big(\rho_s(P^{\ssup{t}})\big)_{s,t\in E}
$$ is regular, where $ P^{\ssup{t}} $ is obtained from $ P^{\ssup{1}} $ by swapping the role of $1$ and $ t $.
\section*{Appendix}
\begin{appendices}

 \section{Pseudo-periodic configurations}\label{pseudoperiodic}
  We define pseudo-periodic configurations as in \cite{DDG12}. We obtain a partition of $ \R^2 $ which is adapted to the Delaunay tessellation. Pick a length scale $ \ell>0 $ and consider the matrix $$ M=\left(\begin{matrix}M_1 & M_2 \end{matrix}\right) =\left(\begin{matrix} \ell & \ell/2\\0 &\sqrt{3}/2 \ell\end{matrix}\right).$$
 Note that $ |M_i|=\ell, i=1,2 $, and $ \angle(M_1,M_2)=\pi/3 $. For each $ (k,l)\in\Z^2 $ we define the cell
 \begin{equation}\label{cell}
 \Delta_{k,l}=\{Mx\in\R^2\colon x-(k,l)\in[-1/2,1/2)^2\}.
 \end{equation} These cells together constitute a periodic partition of $ \R^2 $ into parallelotopes. Let $ B $ be a measurable set of $\O_{\Delta_{0,0}}\setminus\{\emptyset\} $ and
 \begin{equation}\label{ppc}
\Gamma=\{\om\in\O\colon\theta_{Mz}(\om_{\Delta_{k,l}})\in B\mbox{ for all } (k,l)\in\Z^2\} 
 \end{equation} the set of all configurations whose restriction to an arbitrary cell, when shifted back to $ \Delta_{0,0} $, belongs to $ B $. Elements of $ \Gamma $ are called \textbf{pseudo-periodic} configurations. We define marked pseudo-periodic configurations in an analogous way. 
\section{Topology of local convergence}
We write $ \Mcal_1^{\varTheta}(\bO) $ (resp. $ \Mcal_1^{\varTheta}(\O) $) for the set of all shift-invariant probability measures on $ (\bO,\boldsymbol{\Fcal}) $ (resp. $ (\O,\Fcal) $). A measurable function $ f\colon\bO\to\R $ is called local and tame if
$$
f(\bo)=f(\bo_\L)\quad\mbox{ and }\quad |f(\bo)|\le a N_\L(\bo)+b
$$ for all $ \bo\in\bO $ and some $ \L\Subset\R^2 $ and suitable constants $ a,b\ge 0 $. Let $ \Lscr $ be the set of all local and tame functions. The topology of local convergence, or $\Lscr$-topology, on $ \Mcal_1^{\varTheta}(\bO) $ is then defined as the weak$*$ topology induced by $ \Lscr $, i.e., as the smallest topology for which the mappings $ P\mapsto \int f\d P $ with $ f\in\Lscr $ are continuous.
\end{appendices}

\section*{Acknowledgments}
S.A. thanks H.O. Georgii for drawing his attention to the paper \cite{BBD04} and D. Dereudre for helpful discussions and exchange of ideas during a workshop in Leiden. S.A. also thanks the probability group and the department at UBC Vancouver for the warm hospitality during his sabbatical year. M.E. thanks  for the support of his PhD by the EPSRC (Engineering and Physical Sciences Research Council) funded doctoral training centre MASDOC at Warwick.


\begin{thebibliography}{WWW98} 


\bibitem[A15]{A15} \textsc{S.~Adams}, \newblock  A multi-scale approach for geometry-dependent random-cluster percolation, \newblock in preparation (2015).

\smallskip





\bibitem[ACK11]{ACK11}
\textsc{S.~Adams, A.~Collevecchio} and  {\sc W.~K\"onig},
\newblock A variational formula for the free energy of an interacting many-particle system,
\newblock Ann.~Probab. \textbf{39:2}, 683--728 (2011).


\smallskip


\bibitem[BS98]{BS98} \textsc{T.~Benjamini} and \textsc{O.~Schramm}, \newblock Conformal Invariance of Voronoi Percolation, \newblock Commun. Math. Phys. \textbf{197}, 75--107 (1998).


\smallskip

\bibitem[BBD99]{BBD99} \textsc{E.~Bertin, J.M.~Billiot} and \textsc{R.~Drouilhet}, \newblock Existence of Delaunay pairwise Gibbs point process with superstable component, \newblock  Journal of Statistical Physics \textbf{95}, 719--744 (1999).

\smallskip


\bibitem[BBD02]{BBD02} \textsc{E.~Bertin, J.M.~Billiot} and \textsc{R.~Drouilhet}, \newblock Continuum Percolation in The Gabriel Graph, \newblock Advances in Applied Probability, Vol. \textbf{34}, No. 4, 689--701 (2002).





\smallskip


\bibitem[BBD04]{BBD04} \textsc{E.~Bertin, J.M.~Billiot} and \textsc{R.~Drouilhet}, \newblock Phase Transition in the Nearest-Neighbor Continuum Potts Model, \newblock Journal of Statistical Physics \textbf{114}, Nos. 1/2, 79 --100 (2004).



\smallskip


\bibitem[BR06]{BR06}  \textsc{B.~Bollob\'{a}s} and \textsc{O.~Riordan}, \newblock The critical probability for random Voronoi percolation in the plane is $1/2$, \newblock Probab. Theory Relat. Fields \textbf{136}, 417--468 (2006).

\smallskip


\bibitem[BPT]{BPT} \textsc{D.P.~Bourne, M.~Peletier} and \textsc{F.~Theil}, \newblock Optimality of the Triangular Lattice for a Particle System with Wasserstein Interaction, \newblock Commun. Math. Phys. \textbf{329}, 117--140 (2014).



\smallskip 

\bibitem[BKL]{BKL} \textsc{J.~Bricmont, K.~Kuroda} and \textsc{J.L.~Lebowitz}, \newblock The structure of Gibbs states and coexistence for non-symmetric continuum Widom-Rowlinson models, \newblock Prob.~Theory Relat. Field  {\bf 67}, 121--138 (1984).
\smallskip

\bibitem[CCK95]{CCK} \textsc{J.T.~Chayes, L.~Chayes} and \textsc{R.~Koteck\'{y}}, \newblock The Analysis of the Widom-Rowlinson Model by Stochastic Geometric Methods, \newblock Commun. Math. Phys. \textbf{172}, 551--569 (1995).

\smallskip

\bibitem[Der08]{Der08} \textsc{D.~Dereudre}, \newblock Gibbs Delaunay tessellations with geometric hard core conditions, \newblock Journal of Statistical Physics \textbf{131}, 127--151 (2008).


\smallskip




\bibitem[DDG12]{DDG12} \textsc{D.~Dereudre, R.~Drouilhet}  and \textsc{H.O.~Georgii}, \newblock Existence of Gibbsian point processes with geometry-dependent interactions, \newblock Probab. Theory Relat. Fields \textbf{153},  643--670 (2012).

\smallskip

\bibitem[DG09]{DG09} \textsc{D.~Dereudre}  and \textsc{H.O.~Georgii}, \newblock
Variational characterisation of Gibbs measures with Delaunay triangle interaction, \newblock 
Electronic Journal of Probability \textbf{14}, 2438--2462 (2009).

\smallskip

\bibitem[DL11]{DL11} \textsc{D.~Dereudre}  and \textsc{F.~Lavancier}, \newblock Practical Simulation and estimation for Gibbs Delaunay-Voronoi tessellations with geometric hard-core interaction, \newblock Computational Statistics and Data Analysis \textbf{55}, 498--519 (2011).


\smallskip


%

\bibitem[Eye14]{E14} \textsc{M.~Eyers}, \newblock \textit{On Delaunay Random-Cluster Models}, \newblock PhD thesis University of Warwick (2014), available at http://wrap.warwick.ac.uk/67144/.




\smallskip

\bibitem[Geo94]{Geo94}
\newblock {\sc H.-O.~Georgii,}
\newblock Large deviations and the equivalence of ensembles for Gibbssian particle systems with superstable interaction,
\newblock  Prob.~Theory Relat. Fields  {\bf 99}, 171--195 (1994).

\smallskip

\bibitem[Geo88]{G88}
\newblock {\sc H.-O.~Georgii},
\newblock {\it Gibbs Measures and Phase Transitions},
\newblock Berlin: de Gruyter (1988).

\smallskip



\bibitem[GH96]{GH96} \textsc{H.O.~Georgii} and \textsc{O.~H\"aggstr\"om}, \newblock Phase transition in continuum Potts models, \newblock Commun. Math. Phys. \textbf{181}, 507--528 (1996).

\smallskip

\bibitem[GHM]{GHM} \textsc{H.O.~Georgii, O.~H\"aggstr\"om} and \textsc{C.~Maes}, \newblock In: C. Domb and J.L. Lebowitz (eds.) Phase Transitions and Critical Phenomena Vol. \textbf{18}, Academic Press, London, pp. 1--142 ( 2000).




\smallskip


\bibitem[GK97]{GK97} \textsc{H.O.~Georgii} and \textsc{T.~K\"uneth}, \newblock Stochastic comparison of point random fields, \newblock Journal of Applied Probability \textbf{34}, 868--881 (1997).




%
%
%
%
%

\smallskip

\bibitem[Gri94]{G94} \textsc{G.~Grimmett}, \newblock Potts-models and random-cluster processes with many-body interactions, \newblock Journal in Statistical Physics \textbf{75}, 67--121 (1994).




\smallskip

\bibitem[Hag00]{Haggstrom00} \textsc{O.~H\"aggstr\"om}, \newblock Markov random fields and percolation on general graphs, Adv. Appl. Prob. \textbf{32}, 39--66 (2000).

\smallskip

\bibitem[LMP99]{LMP99} \textsc{J.L.~Lebowitz, A.~Mazel} and \textsc{E.~Presutti}, \newblock Liquid-vapor phase transitions for systems with finite range interactions, \newblock Journal in Statistical Physics \textbf{94}, 955--1025 (1999).
\smallskip

\bibitem[LL72]{LL72} \textsc{J.L.~Lebowitz} and \textsc{E.H.~Lieb}, \newblock Phase transition in a continuum classical system with finite interactions, \newblock Phys. Lett. \textbf{39A}, 98--100 (1972).

\smallskip

\bibitem[Lis94]{Lis94} \textsc{D.~Lischinski}, \newblock Incremental Delaunay Triangulation, in \textit{Graphic gems IV}, 47--59, Academic Press (1994).


\smallskip

\bibitem[MSS]{MSS} \textsc{A.~Mazel, Y.~Suhov} and \textsc{I.~Stuhl}, \newblock A classical WR model with q particle types, \newblock arXiv:1311.0020v1 (2013).

\smallskip

\bibitem[MR09]{MR09} \textsc{F.~Merkl} and \textsc{S.W.W.~Rolles}, \newblock Spontaneous breaking of continuous rotational symmetry in two dimensions, \newblock Electron. J. Probab. \textbf{14}, 1705--1726 (2009).



\smallskip

\bibitem[M{\o}94]{M94}
\textsc{J.~M{\o}ller}, \textit{Lectures in Random Voronoi Tessellations}, Springer Lecture Notes in Statistics \textbf{87}, Springer (1994).

\smallskip


\bibitem[T14]{T14} \textsc{V.~Tassion}, \newblock Crossing probabilities for Voronoi percolation,
arXiv :1410.6773v1 (2014).
\smallskip

\bibitem[Rue70]{R70}  \textsc{D.~Ruelle}, \newblock Superstable interactions in classical statistical mechanics, \newblock Commun. Math. Phys. \textbf{18}, 127--159 (1970).

\smallskip


\bibitem[Rue71]{Rue71} \textsc{D.~Ruelle}, \newblock Existence of a phase transition in a continuous classical system, \newblock Phys. Rev.
Lett. \textbf{27}, 1040--1041 (1971).


\smallskip

%
%

\smallskip

\bibitem[WR70]{WR70} \textsc{B.~Widom} and \textsc{J.~Rowlinson}, \newblock New model for the study of liquid-vapour transitions, \newblock 
J. Chem. Phys. \textbf{52}, 1670--1684 (1970).





\end{thebibliography}
\end{document}